\documentclass[10pt]{article}
\usepackage{amsfonts,amsthm,amsmath,amssymb,latexsym,bm,graphics,indentfirst, color}
\newtheorem{theorem}{Theorem}[section]
\newtheorem{lemma}[theorem]{Lemma}
\newtheorem{proposition}[theorem]{Proposition}
\newtheorem{corollary}[theorem]{Corollary}
\newtheorem{remark}[theorem]{Remark}

\newcommand{\D}{\displaystyle}

\newcommand{\p}{\partial}

\newcommand{\reset}{\setcounter{equation}{0}}
\newcommand{\BE}{\begin{equation}}
\newcommand{\BEN}{\begin{equation*}}
\newcommand{\EE}{\end{equation}}
\newcommand{\EEN}{\end{equation*}}
\newcommand{\BL}{\begin{lemma}}
\newcommand{\EL}{\end{lemma}}
\newcommand{\BT}{\begin{theorem}}
\newcommand{\ET}{\end{theorem}}
\newcommand{\BP}{\begin{proposition}}
\newcommand{\EP}{\end{proposition}}
\newcommand{\BC}{\begin{corollary}}
\newcommand{\EC}{\end{corollary}}
\newcommand{\BR}{\begin{remark}}
\newcommand{\ER}{\end{remark}}

\newcommand{\e}{\varepsilon}
\newcommand{\R}{{\mathbb R}}
\newcommand{\C}{{\mathbb C}}
\newcommand{\N}{{\mathbb N}}

\allowdisplaybreaks

\title{Classification and nondegeneracy of $SU(n+1)$ Toda system with singular sources}
\author{Chang-Shou Lin \footnote{Taida Institute of Mathematical of Science, Center of Advanced study in Theoretical Sciences, National Taiwan University, Taipei, Taiwan. {\sl cslin@math.ntu.edu.tw}}\; Juncheng Wei\footnote{Department of Mathematics, Chinese University of Hong Kong, Shatin, Hong Kong.
{\sl wei@math.cuhk.edu.hk}} \; and Dong Ye\footnote{D\'epartement de
Math\'ematiques, UMR 7122, Universit\'e de Metz, B\^{a}t. A, Ile de
Saulcy, 57045 Metz Cedex 1, France. {\sl dong.ye@univ-metz.fr}}}

\begin{document}
\maketitle
\begin{abstract}
We consider the following Toda system
\begin{align*}
\Delta u_i + \D \sum_{j = 1}^n a_{ij}e^{u_j} = 4\pi\gamma_{i}\delta_{0} \;\; \text{in }\mathbb R^2, \quad
\D \int_{\mathbb R^2}e^{u_i} dx < \infty,\;\; \forall\; 1\leq i \leq n,
\end{align*}
where $\gamma_{i} > -1$, $\delta _0$ is Dirac measure at $0$, 
and the coefficients $a_{ij}$ form the standard tri-diagonal
Cartan matrix. 
In this paper, (i)  
we completely classify the solutions and obtain the quantization result:
$$\sum_{j=1}^n a_{ij}\int _{\R^2}e^{u_j} dx = 4\pi (2+\gamma_i+\gamma _{n+1-i}), \;\;\forall\; 1\leq i \leq n.$$
This generalizes the classification result by Jost and Wang for $\gamma_i=0$, $\forall \;1\leq i\leq n$. 
(ii) We prove that if $\gamma_i+\gamma_{i+1}+\ldots+\gamma_j \notin \mathbb Z$ for all $1\leq i\leq j\leq n$, then any solution $u_i$ is \textit{radially symmetric} w.r.t.~$0$. 
(iii) We prove that the linearized equation at any solution is \textit{non-degenerate}. 
These are fundamental results in order to understand the bubbling behavior of the Toda system.
\end{abstract}

\section{Introduction}
In this article, we consider the 2-dimensional (open) Toda system for $SU(n+1)$:
\begin{equation}\label{a4}
\left \{
\begin{array}{l}
\triangle u_i+\displaystyle \sum _{j=1}^na_{ij}e^{u_j}=4\pi \displaystyle \sum _{j=1}^m\gamma _{ij}\delta_{P_j}\quad \mbox{in }\R ^2\\
\D \int_{\R^2}e^{u_i}dx<+\infty
\end{array}\right.
\end{equation}
for $i=1,2,\ldots,n$, where $\gamma _{ij}>-1$, $P_j$ are distinct points and $A=(a_{ij})$ is the Cartan matrix for $SU(n+1)$, given by

\BE\label{cartan} A := (a_{ij}) = \begin{pmatrix} 2 & -1 & 0 & \ldots & 0\\
-1 & 2 & -1 & \ldots & 0\\
0 & -1 & 2 & & 0\\
\vdots& \vdots &&& \vdots\\
0 & \ldots & -1 & 2 & -1\\
 0 & \ldots && -1 & 2
\end{pmatrix} \ . 
\EE
Here $\delta _P$ denotes the Dirac measure at $P$. For $n=1$, system \eqref{a4} is reduced to the Liouville equation 
\BE\label{e103}
\triangle u+e^u=4\pi \displaystyle \sum_{j=1}^m\gamma _j \delta _{P_j}
\EE
which has been extensively studied for the past three decades. The Toda system \eqref{a4} and the Liouville equation \eqref{e103} arise in many physical and geometric problems. For example, in the Chern-Simons theory, the Liouville equation is related to abelian gauge field theory, while the Toda system is related to nonabelian gauge, see \cite{DJPT}, \cite{D}, \cite{GGD}, \cite{L}, \cite{LS}, \cite{M}, \cite{NT}, \cite{non}, \cite{Y1}, \cite{Y} and references therein. On the geometric side, the Liouville equation with or without singular sources is related to the problem of prescribing Gaussian curvature proposed by Nirenberg, or related to the existence of the metrics with conic singularities. 
As for the Toda system, there have been vast articles in the literature to discuss the relationship to holomorphic curves in $\mathbb{CP}^n$, flat $SU(n+1)$ connection, complete integrability and harmonic sequences. For example, see \cite{BJRW}, \cite{BW}, \cite{C}, \cite{CW}, \cite{D1}, \cite{G}, \cite{LS}. In this paper, we want to study the Toda system from the analytic viewpoint. For the past thirty years, the Liouville equation has been extentively studied by the method of nonlinear partial differential equations, see \cite{BM}, \cite{CL}, \cite{CL1}, \cite{CL2}, \cite{Li}, \cite{LY},\cite{NT}, \cite{non},  \cite{pt} and references therein. Recently, the analytic studies of the Toda system can be found in \cite{JLW}, \cite{JW1}, \cite{JW2}, \cite{LL}, \cite{MN}, \cite{non}, \cite{OS}, \cite{WZZ}, \cite{Y1}. For the generalized Liouville system, see \cite{LZ} and \cite{LZ1}.

\medskip
From the pointview of PDE, we are interested not only in the Toda system itself, but also in the case with non-constant coefficients. One of such examples is the Toda system of mean field types:
\BE\label{e103n1}
\Delta u_i(x)+\sum_{j=1}^n a_{ij}\rho_j\left(\frac{h_je^{u_j}}{\D\int_\Sigma h_je^{u_j}}-\frac{1}{|\Sigma|}\right)=4\pi \sum_{j=1}^m\gamma_{ij} \left(\delta_{P_j}-\frac{1}{|\Sigma|}\right),
\EE
where $P_j$ are distinct points, $\gamma_{ij}>-1$ and $h_j$ are positive smooth functions in a compact Riemann surface $\Sigma$. When $n=1$, the equation becomes the following mean field equation:
\BE\label{e104n1}
\Delta u(x)+\rho\left(\frac{he^u}{\D \int_\Sigma he^u}-\frac{1}{\Sigma}\right)=4\pi \sum_{j=1}^m\gamma_j\left(\delta_{P_j}-\frac{1}{|\Sigma|}\right) \quad \mbox{in }\Sigma.
\EE

This type of equations has many applications in different areas of research, and has been extentively investigated.
One of main issues is to determine the set of parameter $\rho$ (non-critical parameters) such that the a priori estimates exist for solutions of equation \eqref{e104n1}. After a priori estimates, we want to compute the topological degree for those non-critical parameters. In this way, we are able to solve the equation \eqref{e104n1} and understand the structure of the solution sets. For the past ten years, those projects have been successfully carried out. See \cite{CL}, \cite{CL1}, \cite{CL2}, \cite{Li}. While carrying out those projects, there often appears a sequence of bubbling solutions and the difficult issue is how to understand the behavior of bubbling solutions near blowup points. For that purpose, the fundamental question is to completely classify all entire solutions of the Toda system with a single singular source:
\BE\label{e103n}
\triangle u_i+\sum_{j = 1}^n a_{ij} e^{u_j}=4\pi \gamma_i\delta _0 \;\; \mbox{in }\R^2,\;\; \D\int_{\R^2}e^{u_i}dx< \infty, \quad 1 \leq i \leq n
\EE
where $\delta _0$ is the Dirac measure at 0, and $\gamma_i > -1$. When all $\gamma_i$ are zero, the classification has been done by Jost-Wang \cite{JW2}. However, when $\gamma_i\neq 0$ for some $i$, the classification has been not solved and has remained a long-standing open problem for many years. It is the purpose of this article to settle this open problem.

\medskip
To state our result, we should introduce some notations.
For any solution $u=(u_1,\cdots,u_n)$ of \eqref{e103n}, we define $U=(U_1,U_2,\cdots,U_n)$ by 
\BE\label{e104n}
U_i=\displaystyle \sum _{j=1}^n a^{ij}u_j
\EE
where $(a^{ij})$ is the inverse matrix of $A$. By \eqref{e104n}, $U$ satisfies
\BE\label{e105n}
\triangle U_i+e^{u_i}=4\pi \alpha_i\delta_0\quad \mbox{in }\R^2,\quad \mbox{where }\;\alpha_i= \sum_{j=1}^n a^{ij}\gamma_j. 
\EE
By direct computations, we have 
\begin{align*}
a^{ij}=\frac{j(n+1-i)}{n+1},\;\;\forall\; n \geq i\geq j\geq 1\quad \mbox{and}\quad 
u_i=\displaystyle \sum_{j=1}^n a_{ij}U_j.
\end{align*}
Our first result is the following classification theorem.

\begin{theorem}\label{thm101n}
Let $\gamma_i>-1$ for $1\leq i \leq n$ and $U=(U_1,\cdots,U_n)$ be defined by \eqref{e104n} via a solution $u$ of \eqref{e103n}. Then $U_1$ can be expressed by 
\BE\label{e107n}
U_1=|z|^{-2\alpha _1}\left(\lambda_0+\displaystyle \sum _{i=1}^n\lambda_i|P_i(z)|^2\right)
\EE
where 
\BE\label{e108n}
P_i(z)=z^{\mu_1+\cdots+\mu_i}+\displaystyle \sum _{j=0}^{i-1}c_{ij}z^{\mu_1+\cdots\mu_j},
\EE
$\mu _i=1+\gamma_i>0$, $c_{ij}$ are complex numbers and $\lambda_i>0$, $0 \leq i\leq n$, satisfy
\BE\label{e109n}
\lambda_0 \cdots \lambda_n=2^{-n(n+1)}\prod_{1\leq i\leq j\leq n}\left(\sum_{k=i}^j\mu_k\right)^{-2} \ .
\EE
Furthermore, if $\mu_{j+1}+\cdots+\mu_{i}\notin\N$ for some $j < i$, then $c_{ij}=0$.
\end{theorem}

In particular, we have the following theorem, generalizing a result by Prajapat-Tarantello \cite{pt} for the singular Liouville equation, $n=1$.
\begin{corollary}\label{cor103}
Suppose $\mu_j+\cdots+\mu_i\notin \N$ for all $1\leq j\leq i \leq n$. Then any solution of \eqref{e103n} is radially symmetric with respect to the origin.
\end{corollary}

We note that once $U_1$ is known, $U_2$ can be determined uniquely by \eqref{e105n}, i.e., $e^{-U_2}=e^{-2U_1}\triangle U_1$. In general, $U_{i+1}$ can be solved via the equation \eqref{e105n} by the induction on $i$. See the formula (\ref{e514}). In the appendix, we shall apply Theorem \ref{thm101n} to give all the explicit solutions in the case of $ n=2$. In section 5, we will prove any expression of \eqref{e107n} satisfying \eqref{e109n} can generate a solution of \eqref{e103n}. See Theorem \ref{lem503}. Thus, the number of free parameters depends on all the Dirac masses $\gamma_j$.
For example if all $\mu_j\in \N$, then the number of free parameters is $n(n+2)$. And if all $\mu_i+\cdots+\mu_j\notin \N$ for $1\leq i\leq j\leq n$, thus the number of free parameters is $n$ only. We let $N(\gamma)$ denote the real dimension of the solution set of the system \eqref{e103n}. 

\medskip
Next, we will show the quantization of the integral of $e^{u_i}$ over $\R^2$ and the non-degeneracy of the linearized system. For the Liouville equation with single singular source:
\begin{align*}
\triangle u+e^u = 4\pi \gamma \delta_0,\quad\int_{\R^2} e^udx <+\infty, \;\;\gamma > -1,
\end{align*}
it was proved in \cite{pt} that any solution $u$ satisfies the following quantization:
\begin{align*}
\int_{\R^2}e^udx=8\pi (1+\gamma ),
\end{align*}
and in \cite{edm} that for any $\gamma \in \N$, the linearized operator around any solution $u$ is nondegenerate.
Both the quantization and the non-degeneracy are important when we come to study the Toda system of mean field type. In particular, this nondegeneracy plays a fundamental role as far as sharp estimates of bubbling solutions of Toda system \eqref{a4}.

\begin{theorem}\label{thm105}
Suppose $u=(u_1,\cdots,u_n)$ is a solution of \eqref{e103n}. Then the followings hold:
\begin{itemize}
\item[(i)] Quantization: we have, for any $1\leq i \leq n$,
$$\sum _{j=1}^na_{ij}\int_{\R^2} e^{u_j}dx =4\pi (2+\gamma_i+\gamma_{n+1-i})$$
and $u_i(z)=-(4+2\gamma_{n+1-i})\log |z|+O(1)$ as $|z|\rightarrow \infty$.
\item [(ii)] Nondegeneracy: The dimension of the null space of the linearized operator at $u$ is equal to $N(\gamma)$.
\end{itemize}
\end{theorem}

In the absence of singular sources, i.e., $\gamma_i=0$ for all $i$, Theorem \ref{thm101n} was obtained by Jost and Wang \cite{JW2}. By applying the holonomy theory, and identifying $S^2=\C\cup\{\infty\}$, they could prove that any solution $u$ can be extended to be a totally unramified holomorphic curve from $S^2$ to $\mathbb{CP}^n$, and then Theorem \ref{thm101n} can be obtained via a classic result in algebraic geometry, which says that any totally unramified holomorphic curve of $S^2$ into $\mathbb{CP}^n$ is a rational normal curve. Our proof does not use the classical result from algebraic geometry. As a consequence, we give a proof of this classic theorem in algebraic geometry by using nonlinear partial differential equations. In fact, our analytic method can be used to prove a generalization of this classic theorem.

\medskip
For a holomorphic curve $f$ of $S^2$ into $\mathbb{CP}^n$, we recall the $k$-th associated curve $f_k:S^2\rightarrow GL(k,n+1)$ for $k=1,2,\cdots,n$ with $f_1=f$ and $f_k=[f\wedge\cdots\wedge f^{(k-1)}]$. A point $p\in S^2$ is called a ramificated point if the pull-back metric $f_k^*(\omega_k)= |z-p| ^{2\gamma _k}h(z) dz\wedge d \bar z$ with $h>0$ at $p$ for some $\gamma_k >0$ where 
\begin{align}
\label{add1}
\mbox{$\omega_k$ is the Fubini-Study metric on $GL(k,n+1)\subseteq \mathbb{CP}^{N_k}$,\ \ $N_k=\begin{pmatrix}
n+1\\k
\end{pmatrix}$}. 
\end{align}
The positive integer $\gamma_k(p)$ is called the ramification index of $f_k$ at $p$. See \cite{GH}.

\begin{corollary}\label{cor106}
Let $f$ be a holomorphic curve of $S^2$ into $\mathbb{CP}^n$. Suppose $f$ has exactly two ramificated points $P_1$ and $P_2$ and $\gamma_j(P_i)$ are the ramification index of $f_j$ at $P_i$, where $f_j$ is the j-th associated curve for $1\leq j\leq n$. Then $\gamma_j(P_1)=\gamma_{n+1-j}(P_2)$. Furthermore, if $f$ and $g$ are two such curves with the same ramificated points and ramification index, then $g$ can be obtained via $f$ by a linear map of $\mathbb{CP}^n$.
\end{corollary}

It is well-known that the Liouville equation as well as the Toda system are completely integrable system, a fact known since Liouville \cite{L1}. Roughly speaking, any solution of \eqref{a4} without singular sources in a simply connected domain $\Omega$ arises from a holomorphic function from $\Omega$ into $\mathbb{CP}^{n}$. See \cite{BJRW}, \cite{BW}, \cite{C}, \cite{CW}, \cite{D1}, \cite{G}, \cite{LS}, \cite{ZS}. For $n=1$,
The classic Liouville theorem says that if a smooth solution $u$ satisfies $\Delta u + e^u = 0$ in a simply connected domain $\Omega \subset \R^2$, then $u(z)$ can be expressed in terms of a holomorphic function $f$ in $\Omega$:
\BE\label{e104}
u(z)=\log \frac{8|f'(z)|^2}{(1+|f(z)|^2)^2}\quad \mbox{in }\; \Omega
\EE
Similarly, system \eqref{a4} has a very close relationship with holomorphic curves in $\mathbb {CP}^n$. Let $F_0$ be a holomorphic curve from $\Omega$ into $\mathbb{CP}^n$. Lift locally $F_0$ to $\C ^{n+1}$ and denote the lift by $\nu=(\nu_0,\nu_1,\ldots,\nu_n)$. The $k$-th associated curve of $F_0$ is defined by
\BE\label{e105}
f_{k}:\Omega \rightarrow G(k,n+1) \subset \mathbb{CP}^{N_k-1},\quad f_{k}(z)= \left[\nu(z)\wedge \nu'(z)\wedge\cdots \nu^{(k-1)}(z)\right],
\EE
where $N_k$ is given by \eqref{add1} and $\nu ^{(j)}$ stand for the $j$-th derivative of $\nu$ w.r.t.~$z$. Let $\Lambda_k=\nu(z)\wedge \cdots \nu^{(k-1)}(z)$. Then the well-known infinitesimal Pl\"{u}cker formulas (see \cite{GH}) is 
\BE\label{e106}
\frac{\partial ^2}{\partial z\partial \bar z}\log \|\Lambda _k\|^2=\frac{\| \Lambda _{k-1}\|^2\| \Lambda _{k+1}\|^2}{\| \Lambda _k\|^4} \quad \mbox{ for } k = 1,2,\cdots ,n,
\EE
where conventionally we put $\| \Lambda _0\| ^2=1$.
Of course, this formula holds only for $\| \Lambda _k\| >0$, i.e.~for all unramificated points. By normalizing $\| \Lambda _{n+1}\| =1$,  letting 
\BE\label{e107}
U_k(z) = -\log \|\Lambda_k(z) \|^2 + k(n-k+1)\log 2,\quad 1\leq k\leq n
\EE
at an unramificated point $z$, and using the fact that $\sum _{1\leq k\leq n} a_{ik} k(n-k+1)=2$, \eqref{e106} gives
\[-\Delta U_i=\exp \left(\displaystyle \sum _{j=1}^na_{ij}U_j\right)\quad \mbox{in }\Omega \setminus\{P_1,\cdots,P_m\}\]
where $\{P_1,\cdots,P_m\}$ are the set of ramificated points of $F_0$ in $\Omega$.
Since $F_0$ is smooth at $P_j$, we have $U_i=-2\alpha _{ij}\log |z-P_j|+O(1)$ near $P_j$. Thus, $U_i$ satisfies
\[\Delta U_i+\exp \left(\sum_{j=1}^na_{ij}U_j\right)=4\pi \sum _{j=1}^n\alpha_{ij} \delta_{P_j}\quad \mbox{in }\Omega .   \]
The constants $\alpha _{ij}$ can be expressed by the total ramification index at $P_j$ by the following arguments.

\medskip
By the Pl\"{u}cker formulas \eqref{e106}, we have
\[ f_i^*(\omega _i)=\frac{\sqrt{-1}}{2}\exp\left(\sum_{j=1}^n a_{ij}U_j\right)dz\wedge d\bar z.\]
Thus, the ramification index $\gamma_{ij}$ at $f_i$ at $P_j$ is 
\BE\label{e115n}
\gamma_{ij}= \sum_{k=1}^n a_{ik}\alpha_{kj}.
\EE
Set
\BE\label{e108}
u_i=\sum _{j=1}^na_{ij}U_j.
\EE
Then it is easy to see that $u_i$ satisfies \eqref{a4} with $\gamma _{ij}$ is the total ramification index of $F_0$ at $P_j$.

\medskip
Conversely, suppose $u=(u_1,\cdots,u_n)$ is a smooth solution of \eqref{a4} in a simply connected domain $\Omega$. We introduce $w_j$ ($0\leq j\leq n$) by 
\BE\label{e109}
u_i=2(w_i-w_{i-1}),\quad \displaystyle \sum _{i=0}^n w_i=0.
\EE
Obviously, $w_i$ can be uniquely determined by $u$ and satisfies
\BE\label{e110}
\left (\begin{array}{c}
w_0\\ \vdots\\w_i\\ \vdots\\ w_n
\end{array}\right )_{z\bar z}=\frac{1}{8}\left( \begin{array}{c}
e^{2(w_1-w_0)}\\\vdots\\ e^{2(w_{i+1}-w_i)}-e^{2(w_i-w_{i-1})}\\ \vdots \\-e^{2(w_n-w_{n-1})}
\end{array}\right). 
\EE
For a solution $(w_i)$, we set 
\begin{align*}
&U=\left(\begin{array}{cccc}
w_{0,z}&0&\ldots&0\\0& w_{1, z}&&0\\
\vdots &&\ddots&\vdots\\0 &0&\ldots& w_{n,z}
\end{array}\right)+\frac{1}{2}
\left( \begin{array}{cccc}
0&0&\ldots&0\\e^{w_1-w_0}&0&&0\\ \vdots&\ddots&\ddots&\vdots\\0&\ldots&e^{w_n-w_{n-1}}&0
\end{array}\right)
\end{align*}
and
\begin{align*}
&V=-\left(\begin{array}{cccc}
w_{0,\bar z}&0&\ldots&0\\0& w_{1, \bar z}&&0\\
\vdots &&\ddots&\vdots\\0 &0&\ldots& w_{n,\bar z}
\end{array}\right)-\frac{1}{2}
\left( \begin{array}{cccc}
0&e^{w_1-w_0}&\ldots&0\\0&0 &&0\\\vdots&\ddots&\ddots&e^{w_n-w_{n-1}}\\0&0 &\ldots&0
\end{array}\right),
\end{align*}
where 
\begin{align*}
w_z =\frac{1}{2}\left(\frac{\partial w}{\partial x}-i\frac{\partial w}{\partial y}\right)\;\; \mbox{and}\;\; w_{\bar z}=\frac{1}{2}\left(\frac{\partial w}{\partial x}+i\frac{\partial w}{\partial y}\right) \quad \mbox{with }\; z=x+iy.
\end{align*}
A straightforward computation shows that $(w_i)$ is a solution of \eqref{e110} if and only if $U$, $V$ satisfy the Lax pair condition: $U_{\bar z}-V_z-[U,V]=0$. Furthermore, this integrability condition implies the existence of a smooth map $\Phi:\Omega \rightarrow SU(n+1,\C)$ satisfying 
\BE\label{e111}\Phi _z=\Phi U,\quad \Phi _{\bar z}=\Phi V\EE
or equivalently, $\Phi $ satisfies $\Phi^{-1}d\Phi =Udz+Vd\bar z$. Let $\Phi=(\Phi_0, \Phi_1\ldots, \Phi_n)$. By \eqref{e111},
\[d\Phi_0= \left(w_{0,z}\Phi _0+\frac{1}{2} e^{w_1-w_0}\Phi_1\right)dz - w_{0,\bar z}\Phi _0d\bar z,\]
which implies
\begin{align}\label{e119}
d(e^{w_0}\Phi _0)&=e^{w_0}d\Phi_0+e^{w_0}\Phi _0dw_0= \left(2w_{0,z}e^{w_0}\Phi _0+\frac{1}{2}e^{w_1}\Phi_1\right)dz.
\end{align}
Therefore, $e^{w_0}\Phi _0$ is a holomorphic function from $\Omega \rightarrow \C ^{n+1}$. We let $\nu(z)=2^{\frac{n}{2}}e^{w_0}\Phi_0$. By using \eqref{e111}, we have $\nu^{(k)}(z)=2^{\frac{n}{2}-k}e^{w_k}\Phi_k$ for $k=1,2,\ldots,n$. Since $w_0+\cdots+w_n=0$, we have $\left\|\nu \wedge \nu '\wedge \cdots \nu^{(n)}(z)\right\|=1$.
Note that 
$$w_0=-\frac{1}{2}\displaystyle \sum_{j=1}^n\frac{(n-j+1)}{n+1}u_j=-\frac{U_1}{2},$$
hence we have  $e^{-U_1}=e^{2w_0}=2^{-n}\| \nu \| ^2$. 
Thus, \eqref{e107} implies $U_1$ is identical to the solution deriving from the holomorphic curve $\nu(z)$. 
Therefore, the space of smooth solutions of the system \eqref{a4} (without singular sources) in a simply connected domain $\Omega$ is identical to the space of unramificated holomorphic curves of $\Omega$ into $\mathbb{CP}^n$.

\medskip
However, if the system \eqref{a4} has singular sources, then $\R^2\setminus\{P_1,\cdots,P_m\}$ is not simply connected. So, it is natural to ask whether in the case $\gamma_{ij}\in \N$, the space of solutions $u$ of \eqref{a4} can be identical to the space of holomorphic curves of $\R^2$ into $\mathbb{CP}^n$ which ramificates at $P_1,\cdots,P_m$, with the given ramification index $\gamma_{ij}$ at $P_j$. The following theorem answers this question affirmatively.
\begin{theorem}\label{thm101}
Let $\gamma_{ij} \in \mathbb{N}$ and $P_j\in \mathbb{R}^2$. Then for any solution $u$ of \eqref{a4}, there exists a holomorphic curve $F_0$ of $\C$ into $\mathbb{CP}^n$ with ramificated points $P_j$ and the total ramification index $\gamma _{ij}$ at $P_j$ such that for $1 \leq k \leq n$,
\[e^{-U_k}=2^{-k(n+1-k)}\left\|\nu(z)\wedge \cdots \wedge \nu^{(k-1)}(z)\right\| ^2\quad \mbox{in }\; \C \setminus\{P_1,\cdots,P_m\}\]
where $\nu(z)$ is a lift of $F_0$ in $\C^{n+1}$ satisfying
\[\left\| \nu(z)\wedge \cdots \wedge \nu ^{(n)}(z)\right\| =1.\]
Furthermore, $F_0$ can be extended smoothly to a holomorphic curve of $S^2$ into $\mathbb{CP}^n$.
\end{theorem}

We note that if equation \eqref{a4} is defined in a Riemann surface rather than $\C$ or $S^2$, then the identity of the solution space of \eqref{a4} with holomorphic curves in $\mathbb{CP}^n$ generally does not hold. For example, if the equation \eqref{a4} is defined on a torus, then even for $n=1$, a solution of \eqref{a4} would be not necessarily associated with a holomorphic curve from the torus into $\mathbb{CP}^n$. See \cite{LW}.

\medskip
The paper is organized as follows. In section 2, we will show some invariants associated with a solution of the Toda system. Those invariants allows us to classify all the solutions of \eqref{e103n} without singular sources, thus it gives another proof of the classification due to Jost and Wang. Those invariants in section 5 can be extended to be meromorphic invariants for the case with singular sources. By using those invariants, we can prove $e^{-U_1}$ satisfies an ODE in $\C^*:= \C\setminus\{0\}$, the proof will be given in section 5. In section 4 and section 6, we will prove the quantization and the non-degeneracy of the linearized equation of \eqref{e103n} for the case without or with singular sources. In the final section, we give a proof of Theorem \ref{thm101}. Explicits solutions in the case of $SU(3)$ are given in the appendix.

\bigskip
\noindent
{\bf Acknowledgments:} The research of J.W. is partially
supported by a research Grant from GRF of Hong Kong  and a Joint
HK/France Research Grant. D.Y. is partly supported by the French ANR
project referenced ANR-08-BLAN-0335-01, he would like to thank
department of mathematics of CUHK for its hospitality.

\section{Invariants for solutions of Toda system}
\reset
In this section, we derive some invariants for the Toda system. Denote
$A^{-1} = (a^{jk})$, the inverse matrix of $A$. Let \BE\label{a7} U_j = \sum_{k = 1}^n a^{jk}u_k, \quad \forall\;
1\leq j \leq n.\EE  Since $\Delta=4 \partial_{z\bar z}$, it is easy to
see that the system (\ref{e103n}) is equivalent to for all $1\leq i \leq n$,
\begin{align*}
-4U_{i, z\bar z} = {\rm exp}\left({\sum_{j = 1}^n a_{ij}U_j}\right) - 4\pi\alpha_i\delta_0 \;\; \mbox{in }\R^2,
\quad \int_{\R^2}  {\rm exp}\left({\sum_{j = 1}^n a_{ij}U_j}\right)dx
< \infty.\end{align*}
where $\alpha_i = \sum_{1\leq j \leq n} a^{ij}\gamma_j$ for $1\leq i \leq n$.
Define \BE
\label{W} 
W_1^j =- e^{U_1}\left(e^{-U_1}\right)^{(j+1)}\;\; \text{for }1\leq j\leq n \quad \text{and}\quad 
W_{k+1}^j =-\frac{W_{k, \bar z}^j}{U_{k, z\bar z}} \;\; \mbox{ for $1\leq k \leq j-1$}.
\EE
We will prove that all these quantities $W^j_k$, $1\leq k\leq j \leq n$, are invariants for solutions of $SU(n+1)$, more precisely, $W^j_k$ are a part of some specific holomorphic or meromorphic functions, which are determined explicitly by the Toda system. 

\BL\label{lem21}
For any classical solution of \eqref{a4}, there holds:
\BE\label{e204}
W_k^k=\sum_{i = 1}^k (U_{i, zz} -U_{i,z}^2)+\sum_{i = 1}^{k-1} U_{i, z}U_{i+1, z}\quad \text{for } 1\leq k\leq n,
\EE
\BE\label{e205}
W^k_{k, \bar z} =-U_{k,z\bar{z}}U_{k+1,z} \quad \text{ for } 1\leq k\leq n-1,\\[10pt]
\EE
\BE\label{e206}
W_k^j=(U_{k-1, z} -U_{k,z})W_k^{j-1}+W_{k,z}^{j-1}+W_{k-1}^{j-1} \quad \text{ for } 1\leq k < j \leq n.
\EE
where for convenience $U_0 = 0$ and $W_0^j = 0$ for all $j$. 
\EL

\begin{proof}
First, we show that \eqref{e204} implies \eqref{e205}.
By the equation for $U_j$,
\BE\label{e208}
U_{j,z\bar{z}z}=U_{j,z\bar{z}}(2U_{j,z}-U_{j+1,z}-U_{j-1,z}), \quad \forall \; 1 \leq j \leq n, 
\EE
where for the convenience, $U_{n+1}=0$ is also used. Thus, 
\begin{align}
\label{e207}
\begin{split}
-U_{j,z\bar{z}}U_{j+1,z} + U_{j-1, z\bar z}U_{j, z} & = U_{j,z\bar{z}z}- U_{j,z\bar{z}}\left(2U_{j,z}-U_{j-1,z}\right) + U_{j-1, z\bar z}U_{j, z}\\
& = \left(U_{j,zz}-U_{j,z}^2 + U_{j,z}U_{j-1,z}\right)_{\bar z}.
\end{split}
\end{align}
Taking the sum of \eqref{e207} for $j$ from $1$ to $k$, we get
\[-U_{k,z\bar{z}}U_{k+1,z}= \sum_{j=1}^k \left(U_{j,zz}-U_{j,z}^2+ U_{j,z}U_{j-1,z}\right)_{\bar{z}} = W_{k, \bar z}^k\]
where \eqref{e204} is used.

\medskip
Next, we will prove \eqref{e204}-\eqref{e206} by the induction on $k$. Obviously, \eqref{e204} holds for $k=1$. By the definition of $W_1^j$, for $j\geq 2$, we have
\[
W_1^j=-e^{U_1}(e^{-U_1})^{(j+1)}=e^{U_1}\left(e^{-U_1}W_1^{j-1}\right)_z=W_{1,z}^{j-1}-W_1^{j-1}U_{1,z},
\]
which is \eqref{e206} for $k=1$. To compute $W_{k+1}^{k+1}$, \eqref{e206} with index $k$ implies 
\[
-U_{k,z\bar{z}}W_{k+1}^{k+1}=W_{k,\bar{z}}^{k+1}=(U_{k-1,z\bar{z}}-U_{k,z\bar{z}})W_k^k+(U_{k-1,z}-U_{k,z})W_{k,\bar{z}}^k+W_{k,z\bar{z}}^k+W_{k-1,\bar{z}}^k,
\]
Since $U_{k-1,z\bar{z}}W_k^k+W_{k-1,\bar{z}}^k=0$, the above identity leads by \eqref{e205} with index $k$,
\begin{align*}
W_{k,\bar{z}}^{k+1}= & \;-U_{k,z\bar{z}}W_k^k-(U_{k-1,z}-U_{k,z})U_{k,z\bar{z}}U_{k+1,z} -(U_{k,z\bar{z}}U_{k+1,z})_z\\
= & \;-U_{k,z\bar{z}}W_k^k-(U_{k-1,z}-U_{k,z})U_{k,z\bar{z}}U_{k+1,z} - U_{k,z\bar{z}}(2U_{k,z}-U_{k+1,z}-U_{k-1,z})U_{k+1,z}\\
&\; -U_{k,z\bar{z}}U_{k+1,zz}\\
= & \;-U_{k,z\bar{z}}\left(W_k^k+U_{k+1,zz}-U_{k+1,z}^2+U_{k+1,z}U_{k,z}\right)
\end{align*}
where \eqref{e208} is used. Hence
\[W_{k+1}^{k+1}=W_k^k+U_{k+1,zz}-U_{k+1,z}^2+U_{k+1,z}U_{k,z},\]
and then \eqref{e204} is proved for $k+1$.

\medskip
To compute $W_{k+1}^j$ for $j\geq k+2$, we have $j-1\geq k+1$ and by similar calculations:
\begin{align*}
W_{k,\bar{z}}^j = & \;(U_{k-1,z\bar{z}}-U_{k,z\bar{z}})W_k^{j-1}+(U_{k-1,z}-U_{k,z})W_{k,\bar{z}}^{j-1}+W_{k,z\bar{z}}^{j-1}+W_{k-1,\bar{z}}^{j-1}\\
= & \; -U_{k,z\bar{z}}W_k^{j-1}-(U_{k-1,z}-U_{k,z})U_{k,z\bar{z}}W_{k+1}^{j-1}- \left(U_{k,z\bar{z}}W_{k+1}^{j-1}\right)_z\\
= & \; -U_{k,z\bar{z}}W_k^{j-1}-(U_{k-1,z}-U_{k,z})U_{k,z\bar{z}}W_{k+1}^{j-1}-U_{k,z\bar{z}}(2U_{k,\bar{z}}-U_{k+1,z}-U_{k-1,z})W_k^{j-1}\\
& \;-U_{k,z\bar{z}}W_{k+1,z}^{j-1}\\
= & \; -U_{k,z\bar{z}}\left[(U_{k,z}-U_{k+1,z})W_{k+1}^{j-1}+W_{k+1,z}^{j-1}W_k^{j-1}\right],
\end{align*}
which leads to
\[W_{k+1}^j=(U_{k,z}-U_{k+1,z})W_{k+1}^{j-1}+W_{k+1,z}^{j-1}+W_k^{j-1}\]
Therefore, Lemma \ref{lem21} is proved.
\end{proof}

\section{Classification of solutions of $SU(n+1)$ with $m = 0$}
\reset
Here we show a new proof of the classification
result of Jost-Wang \cite{JW2}. That is, all classical solutions of
\eqref{a4} with $m=0$ is given by a $n(n+2)$ manifold $\cal M$. Our idea is to
use the invariants $W_j^n$ for solutions of $SU(n+1)$. Consider
\begin{align}
\label{todaws}
- \Delta u_i = \D \sum_{j = 1}^n a_{ij}e^{u_j} \;\; \text{in }\mathbb R^2, \quad
\D \int_{\mathbb R^2}e^{u_i}dx <\infty, \quad \forall\; 1\leq i \leq n.
\end{align}
\BT \label{n11} 
For any classical solution of \eqref{todaws}, let $U_j$,
$W_j^n$ be defined by \eqref{a7} and \eqref{W}, then \begin{align*}
W_j^n \equiv 0 \quad \mbox{in } \R^2, \qquad \forall\; 1 \leq j \leq
n.\end{align*} \ET \BR The fact $W_n^n = 0$ has been proved by Jost
and Wang in an equivalent form, which is just the function $f$ in
the proof of Proposition 2.2 in \cite{JW2}. \ER

\begin{proof}
The proof is based on the following observation: \BE\label{n4}
W^n_{n, \bar z} = 0 \;\; \mbox{in $\R^2$} \quad \mbox{for any solution of \eqref{todaws}.} \EE In fact, using formula \eqref{e204} and the equations of $U_i$,
\begin{align}
\label{add9}
\begin{split}
W_{n, \bar z}^n & = \sum_{i = 1}^n (U_{i, z\bar z})_z - 2\sum_{i =
1}^n U_{i, z}U_{i, z\bar z} + \sum_{i = 1}^{n-1} \left(U_{i, z\bar
z}U_{i+1, z} + U_{i, z}U_{i+1, z\bar z}\right)\\
& = \sum_{i = 1}^n U_{i, z\bar z}\left[\sum_{j = 1}^n
\left(a_{ij}U_{j, z}\right) - 2U_{i, z} + U_{i+1,
z} + U_{i-1, z}\right]\\
& = 0.
\end{split}
\end{align}
Here we used again the convention $U_0 = U_{n+1} = 0$ for $SU(n+1)$.

\medskip
Furthermore, $e^{u_i} \in L^1(\R^2)$ implies that for any $\epsilon >0$, there exists $R_{\epsilon}>0$ such that
\[\int_{\mathbb R^2\setminus B_{R_{\epsilon}}} e^{u_i}dz \leq \epsilon ,\quad 1\leq i \leq n\]
For sufficient small $\epsilon >0$, applying Brezis-Merle's argument \cite{BM} to the system $u_i$, we can prove $u_i(z)\leq C$ for $|z|\geq R_{\epsilon }$, i.e.~$u_i$ is bounded from the above over $\C$.
Thus, $u_i$ can be represented by the following integral formulas:
\BE\label{e303}
u_i(z)=\frac{1}{2\pi }\int _{\R^2}\log \frac{|z'|}{|z-z'|}\displaystyle \sum ^n_{j=1} a_{ij}e^{u_j(z')}dz'+c_i, \quad\forall\; 1\leq i\leq n,
\EE
for some real constants $c_i$.

\medskip
This gives us the asymptotic
behavior of $u_i$ and their derivatives at infinity. In particular,
for any $k \geq 1$, $\nabla^k u_i = O\left(|z|^{-k}\right)$ as $|z|$ goes to
$\infty$. So $\nabla^k U_i =
O\left(|z|^{-k}\right)$ as $|z| \to \infty$, for $k \geq 1$.
Therefore, $W_n^n$ is a entire holomorphic function, which tends to
zero at infinity, so $W_n^n \equiv 0$ in $\R^2$ by classical
Liouville theorem. As $W_{n-1, \bar z}^n = -U_{n-1, z\bar z}W_n^n$,
we obtain $W_{n-1, \bar z}^n = 0$ in $\R^2$. By \eqref{e204} and \eqref{e206}, it is not difficult to
see that for $1\leq i \leq n-1$, $W_i^n$ are also polynomials of
$\nabla^k U_i$ with $k \geq 1$, so they tend to $0$ at infinity,
hence $W_{n-1}^n = 0$ in $\R^2$. We can complete the proof of
Theorem \ref{n11} by induction. \end{proof}
Futhermore, we know that $e^{-U_1}$ can be computed as a square of some holomorphic curves in $\mathbb{CP}^n$, see the Introduction. Thus, there is a holomorphic map $\nu(z)=(\nu_0(z),\ldots,\nu_n(z))$ from $\C$ into $\C ^{n+1}$ satisfying 
\[\left\|\nu \wedge \nu' \cdots \wedge \nu ^{(n)}(z)\right\|=1 \quad \mbox{and}\quad e^{-U_1(z)}=\sum _{i=0}^n|\nu _i(z)|^2\;\; \mbox{in }\;  \C.\]
Since $W_1^n \equiv 0$ in $\R^2$ yields $(e^{-U_1})^{(n+1)}=0$, we have $\nu _i^{(n+1)}(z)=0$. By the asymptotic behavior of $u_i$, we know that $e^{-U_1}$ is of polynomial growth as $|z|\rightarrow \infty$. Hence $\nu _i(z)$ is a polynomial and $\nu_0,\cdots, \nu_n$ is a set of fundamental holomorphic solutions of $f^{(n+1)}=0$. Thus
\begin{align}
\label{add2}
\nu_i(z)= \sum_{j = 0}^n c_{ij}z^j\quad \mbox{with }\;\; {\rm det}(c_{ij})\neq 0.
\end{align}
By a linear transformation, we have 
\[\nu(z)=\lambda(1,z,z^2,\cdots,z^n), \quad \lambda \in \C\]
and $[\nu ]$ is the rational normal curve of $S^2$ into $\mathbb{CP}^n$.
Hence we have proved the classification theorem of Jost and Wang.

\BR
Here we use the integrability of the Toda system. In section 5, we actually prove the classification theorem without use of the integrability.
\ER

\BR The invariants $W_j^n$ are called $W$-symmetries or conservation
laws, see \cite{LS}. It is claimed that for the Cartan matrix there
are $n$ linearly independent $W$-symmetries, see \cite{ZS}. However,
as far as we are aware, we cannot find the explicit formulas in the
literature (except for $n=2$ \cite{WZZ}). Here we give explicit formula for the
$n$ invariants. \ER

\section{Nondegeneracy of solutions of $SU(n+1)$ without sources}
\reset
Let $\mathcal M$ be the collection of entire solution of \eqref{todaws}. In the previous section, we know that $\mathcal M$ is a smooth manifold of $n(n+2)$ dimension. Fixing a solution $u=(u_1,\cdots,u_n)$ of \eqref{todaws}, we consider $LSU(n+1)$, the linearized system of \eqref{todaws} at $u$:
\BE\label{e401}
\triangle \phi _i+\sum_{j=1}^n a_{ij} e^{u_j}\phi_j=0\quad \mbox{in }\R^2.
\EE
Let $s\in \R$ be any parameter appearing in \eqref{add2} and $u(z;s)$ be a solution of \eqref{todaws} continuously depending on $s$ such that $u(z;0)=u(z)$. Thus $\phi(z)=\frac{\partial}{\partial s}u(z;s)|_{s=0}$ is a solution of \eqref{e401} satisfying $\phi\in L^{\infty}(\R^2)$. Let $T_u\mathcal M$ denote the tangent space of $\mathcal M$ at $u$. The nondegeneracy of the linearized system is equivalent to showing that any bounded solution $\phi =(\phi_1,\cdots,\phi_n)$ of \eqref{e401} belongs to this space.

\begin{theorem}\label{thm401}
Suppose $u$ is a solution of \eqref{todaws} and $\phi$ is a bounded solution of \eqref{e401}. Then $\phi \in T_u\mathcal M$.
\end{theorem}

\begin{proof}
For any solution $\phi =(\phi _1,\ldots,\phi_n)$ of \eqref{e401}, we define
\BE\label{b2} \Phi_j = \sum_{k = 1}^n a^{jk}\phi_k, \quad
\forall\; 1\leq j \leq n.\EE 
We have readily that bounded $(\phi_i)$ solves
\eqref{e401} if and only if $(\Phi_i)$ is a solution of \BE\label{b3}
-4\Phi_{i, z\bar z} = {\rm exp}\left({\sum_{j = 1}^n
a_{ij}U_j}\right)\times \sum_{j = 1}^n  a_{ij}\Phi_j\;\;
\text{in }\mathbb R^2, \quad \Phi_i \in L^\infty(\R^2)\quad \forall\; 1\leq i \leq
n. \EE 
Our idea is also to find some invariants which characterize all
solutions of \eqref{b3}. Indeed, we find them by linearizing the
above quantities $W_k^n$ for $U_i$. Let \begin{align*} Y_1^n =
e^{U_1}\left[\left(e^{-U_1}\Phi_1\right)^{(n+1)} -
\left(e^{-U_1}\right)^{(n+1)}\Phi_1\right] \end{align*} and
\begin{align*}
Y_{k+1}^n = - \frac{Y_{k, \bar z}^n + W_{k+1}^n\Phi_{k, z\bar
z}}{U_{k, z\bar z}}\quad \mbox{for $1\leq k \leq n-1$}.\end{align*}
The quantities $Y_k^n$ are well defined and we can prove by
induction the following formula: {\it With any solutions of
$LSU(n+1)$}, there hold 
\begin{align*} & Y_1^{n} = Y_{1, z}^{n-1} -
Y_1^{n-1}U_{1, z} - W_1^{n-1}\Phi_{1, z}\\ & Y_k^{n} =
\left(U_{k-1, z} - U_{k, z}\right) Y_k^{n-1} + Y_{k, z}^{n-1} + Y_{k - 1}^{n-1}
+ \left(\Phi_{k-1, z} - \Phi_{k, z}\right) W_k^{n-1},\quad \mbox{for $2
\leq k \leq n$}.
\end{align*}
Moreover, for any solution of
\eqref{b3}, we have
\begin{align}
\label{n7} Y_n^n = \sum_{i = 1}^n \Phi_{i, zz} - 2\sum_{i = 1}^n
U_{i, z}\Phi_{i, z} + \sum_{i = 1}^{n-1} \left(\Phi_{i, z}U_{i+1, z}
+ U_{i, z}\Phi_{i+1, z}\right).
\end{align}
The proof is very similar as above for $W_j^n$, since each quantity $Y_j^n$ is
just the {\it linearized} version of $W_j^n$ with respect to
$(U_i)$, as well as the involved equations, so we leave the details
for interested readers. 

\medskip
Applying the equations \eqref{b3}, it can be checked
easily that
\begin{align*} \mbox{$Y_{n, \bar z}^n  = 0$\;\; in $\R^2$, \;\; for any
solution of $LSU(n+1)$ \eqref{e401}.}
\end{align*}
Using the classification of $u_i$ in section 3 (see also \cite{JW2}), we know that $e^{u_i} =
O\left(z^{-4}\right)$ at $\infty$. Since $\phi_i \in L^\infty(\R^2)$,  the function $ \sum_{1\leq j \leq n}
a_{ij} e^{u_j} \phi_j \in L^1 (\R^2)$. As before, we can express
$\phi_i$ by integral representation and prove that
$\lim_{|z|\to\infty}\nabla^k\phi_i = 0$ for any $k \geq 1$. Hence
$\lim_{|z|\to\infty}\nabla^k\Phi_i = 0$ for any $k \geq 1$.

\medskip
By similar argument as above, this implies that $Y_n^n = 0$ in $\R^2$ for any
solution of \eqref{b3}, and we get successively $Y_k^n = 0$ in
$\R^2$ for $1\leq k \leq n-1$, recalling just $Y_{k, \bar z}^n =
-U_{k, z\bar z}^nY_{k+1}^n - \Phi_{k, z\bar z}W_{k+1}^n$ and $W_j^n
= 0$ in $\R^2$ for any classical solution of \eqref{todaws}. Since
\begin{align*}
0 = Y_1^n = e^{U_1}\left(e^{-U_1}\Phi_1\right)^{(n+1)} + W_1^n\Phi_1
= e^{U_1}\left(e^{-U_1}\Phi_1\right)^{(n+1)},
\end{align*}
we conclude then $\left(e^{-U_1}\Phi_1\right)^{(n+1)} = 0$ in
$\R^2$. By the growth of real function $e^{-U_1}\Phi_1$, we get
\begin{align*}
e^{-U_1}\Phi_1 = \sum_{i, j = 0}^n b_{ij}z^i\bar z^j
\end{align*}
with $b_{ij} = \overline{b_{ji}}$ for all $0\leq i, j\leq
n$. This yields \begin{align*} \Phi_1 \in \mathcal L = \left\{
e^{U_1}\left[\sum_{i, j = 0}^n b_{ij} z^i\bar z^j\right], \;
b_{ij} \in \C, \; b_{ij} = \overline{b_{ji}},\;
\forall\; 0\leq i, j\leq n \right\},\end{align*} a linear space of
dimension $(n+1)^2$. Once $\Phi_1$ is fixed, as $-\Delta\Phi_1 =
e^{u_1}(2\Phi_1 - \Phi_2)$ in $\R^2$, $\Phi_2$ is uniquely
determined, successively all $\Phi_i$ are uniquely determined, so is
$\phi_i$.

\medskip
Moreover, the expression of $e^{-U_1}$ given by the last section yields
that the constant functions belong to $\mathcal L$. If $\Phi_1 \equiv
\ell_1 \in \R$, by equations \eqref{b3}, successively we obtain
$\Phi_i \equiv \ell_i \in \R$ for all $2\leq i \leq n$. Using again
the system \eqref{b3}, we must have
$$\sum_{j = 1}^n a_{ij}\ell_j = 0,\quad \forall\; 1 \leq i \leq n,$$
which implies $\ell_j = 0$ for any $1\leq j \leq n$, hence $(\Phi_i)$
can only be the trivial solution. Therefore, we need only to
consider $\Phi_1$ belonging to the algebraic complementary of $\R$
in $\mathcal L$, a linear subspace of dimension $n(n+2)$.

\medskip
Finally, it is known that $T_u\mathcal M$, the tangent space of $u = (u_i)$ to
the solution manifold $\cal M$ provides us a $n(n+2)$ dimensional
family of bounded solutions to $LSU(n+1)$, so we can conclude that all the
solutions of \eqref{e401} form exactly a linear space of dimension
$n(n+2)$. Theorem \ref{thm401} is then proved. \end{proof}

\BR
We can remark by the proof that Theorem \ref{thm401} remains valid if we relax the condition $\phi_i \in L^\infty(\R^2)$ to the growth condtion $\phi_i(z) = O(|z|^{1 +\alpha})$ at infinity with $\alpha \in (0, 1)$.
\ER

\section{Classification of singular Toda system with one source}
\reset
For the Toda system $SU(n+1)$ with one singular source \eqref{e103n}, denote 
$A^{-1} = (a^{jk})$, the inverse matrix of $A$ and define as before
\begin{align}
\label{singSU}
U_j = \sum_{k = 1}^n a^{jk}u_k, \quad \alpha_j = \sum_{k = 1}^n a^{jk}\gamma_k\quad \forall\;
1\leq j \leq n.
\end{align}
where $u=(u_1,...,u_n)$ is a solution of \eqref{e103n}. So
\begin{align}\label{e502}
-\Delta U_i = {\rm
exp}\left({\sum_{j = 1}^n a_{ij}U_j}\right) - 4\pi\alpha_i\delta_0
\end{align}
with
\begin{align*}
 \int_{\R^2} {\rm exp}\left({\sum_{j = 1}^n
a_{ij}U_j}\right)dx = \int_{\R^2} e^{u_i}dx < \infty, \;\; \forall\; i.
\end{align*}

In this section, we will completely classify all the solutions of equation \eqref{e103n}, and prove in the next section the nondegenerency of the corresponding linearized system. Here is the classification result.
\begin{theorem}
\label{mainsingular}
Suppose that $\gamma_i >-1$ for $1\leq i\leq n$, and $U=(U_1,...,U_n)$ is a solution of \eqref{e502}, then we have 
\begin{align}
\label{equaU1}
|z|^{2\alpha_1}e^{-U_1} = \lambda_0 + \sum_{1 \leq i \leq n} \lambda_i |P_i(z)|^2  \quad\mbox{in }\C^*
\end{align}
where 
\begin{align}
\label{pi}
\lambda_i \in \R, \quad P_i(z) = c_{i0}+ \sum_{j = 1}^{i-1} c_{ij}z^{\mu_1+\mu_2+\ldots+\mu_j} + z^{\mu_1+\mu_2+\ldots+\mu_i}, \quad  c_{ij} \in \C.
\end{align}
Moreover, $\lambda_i$ verifies the following necessary and sufficient conditions: 
\begin{align}
\label{equal}
\lambda_i > 0, \quad \lambda_0\lambda_1\cdots\lambda_n = 2^{-n(n+1)} \times \prod_{1 \leq i\leq j\leq n}\left(\sum_{k= i}^j\mu_k\right)^{-2}.
\end{align}
Conversely, $U_1$ defined by \eqref{equaU1}-\eqref{equal} generates a solution $(U_i)$ of \eqref{e502}.
\end{theorem}

\medskip
The proof of Theorem \ref{mainsingular} is divided in several steps. Suppose $U=(U_1,...,U_n)$ is a solution of \eqref{e502}.
\subsection{Step 1}
We will prove that $e^{-U_1} = f$ verifies the differential equation as follows:
\begin{align}
\label{equaf}
f^{(n+1)} + \sum_{k = 0}^{n-1} \frac{w_k}{z^{n+1 - k}}f^{(k)} = 0\quad \mbox{in } \C^*,\end{align}
where $w_k $ are real constants only depending on all $\gamma _i$ and $f^{(i)}$ denotes the $i$-th order derivative of $f$ w.r.t.~$z$. 

\begin{lemma}\label{lem502n}
Let $(U_j)$ be given by \eqref{mainsingular}, with $(u_i)$ a solution of \eqref{e103n}. Define $Z_n = W_n^n$ and by iteration
\begin{align}\label{e507}
Z_k = W_k^n + U_{k, z}Z_{k +1} + \sum_{j = k}^{n-2}W_k^j Z_{j+2}, \quad \forall\;  k = n-1, n-2, \ldots, 1.
\end{align}
Then $Z_k$ are holomorphic in $\C^*$. More precisely, there exist $w_k\in \C$ such that 
$$Z_k =  \frac{w_k}{z^{n+2 - k}}\;\; \mbox{ in }\; \C^*, \quad \mbox{for any }\; 1\leq k \leq n,$$
where $w_k$ only depends on $\gamma _j$.
\end{lemma}
Here $W_k^j$ ($1\leq k \leq j\leq n$), considered as functional of $(U_1, U_2, \ldots U_n)$ and their derivatives, are the invariants constructed in section 2 for Toda system $SU(n+1)$.

\medskip
\noindent
\begin{proof}
 First, we recall that
\begin{align}\label{e508} W_1^m =
- e^{U_1}\left(e^{-U_1}\right)^{(m+1)}\;\; \mbox{for $1 \leq m \leq n$}, \quad W_{k+1}^m =
-\frac{W_{k, \bar z}^m}{U_{k, z\bar z}} \;\; \mbox{ for $1\leq k \leq
m-1$}.\end{align}
Using \eqref{add9}, $Z_n$ is holomorphic in $\C^*$ and by Lemma \ref{lem21}
\begin{align*}
W_{k, \bar z}^k = -U_{k, z\bar z}U_{k+1, z}, \quad \mbox{for any }\; 1 \leq k \leq n-1.
\end{align*}
Consequently,  in $\C^*$ there holds by \eqref{e508},
$$0 = W_{n-1, \bar z}^n  + U_{n-1, z\bar {z}}W_n = W_{n-1, \bar {z}}^n  + U_{n-1, z\bar z}Z_n = \left (W_{n-1}^n  + U_{n-1, z}Z_n \right)_{\bar z} = Z_{n-1,\bar{z}},$$
So $Z_{n-1}$ is also holomorphic in $\C^*$. Suppose that $Z_{\ell+1}$ are holomorphic in $\C^*$ for $k \leq \ell \leq n - 2$, then we have in $\C^*$,
\begin{align*}
Z_{k, \bar z} & = \left(W_k^n + U_{k, z}Z_{k +1} + \sum_{j = k}^{n-2}W_k^j Z_{j+2}\right)_{\bar z}\\
& = W_{k, \bar z}^n + U_{k, z\bar z}Z_{k +1} + W_{k, \bar z}^kZ_{k+2} + \sum_{j = k+1}^{n-2}W_{k, \bar z}^j Z_{j+2}\\
& = - U_{k, z\bar z}W_{k+1}^n + U_{k, z\bar z}Z_{k +1} -U_{k, z\bar z}U_{k+1, z}Z_{k+2} - \sum_{j = k+1}^{n-2}U_{k, z\bar z}W_{k+1}^jZ_{j+2}\\
& =  U_{k, z\bar z}\left(Z_{k +1} - W_{k+1}^n - U_{k+1, z}Z_{k+2} - \sum_{j = k+1}^{n-2}W_{k+1}^jZ_{j+2}\right) = 0.
\end{align*}
The last line comes from the definition of $Z_{k+1}$. Thus, $Z_k$ is holomorphic in $\C ^*$ for all $1 \leq k \leq n$.

\medskip
Next, we want to show that 
\BE\label{e509}
Z_k=\frac{w_k}{z^{n+2-k}}
\EE
for some real constant $w_k$ depending on $\gamma _j$. Define
\begin{align}
\label{add3}
V_j = U_j - 2\alpha_j\log |z|, \quad \forall\;
1\leq j \leq n.
\end{align}
So
\begin{align*}
-\Delta V_i = -4U_{i, z\bar z} + 4\pi\alpha_i\delta_0 & = {\rm
exp}\left({\sum_{j = 1}^n a_{ij}U_j}\right) + 4\pi\alpha_i\delta_0 -
4\pi\sum_{j = 1}^n \left(a^{ij}\gamma_j\delta_0\right)\\
& =
|z|^{2\gamma_j}{\rm exp}\left({\sum_{j = 1}^n a_{ij}V_j}\right)
\end{align*}
with
\begin{align*}
\int_{\R^2} |z|^{2\gamma_j}{\rm exp}\left({\sum_{j = 1}^n
a_{ij}V_j}\right) dx = \int_{\R^2} {\rm exp}\left({\sum_{j = 1}^n
a_{ij}U_j}\right) dx = \int_{\R^2} e^{u_i}dx < \infty, \;\; \forall\;1\leq  i\leq n.
\end{align*}
As $\gamma_i > -1$, applying Brezis-Merle's argument in \cite{BM} to the system of $V_i$, we have $V_i \in C^{0, \alpha}$ in $\C$ for some $\alpha \in (0, 1)$ and they are upper
bounded over $\C$. This implies that we can express $V_i$ by the
integral representation formula. Moreover, by scaling argument and elliptic estimates, we have for all $1 \leq i \leq n$,
\begin{align}
\label{add4}
\nabla^kV_i(z) = O\left(1 + |z|^{2+2\gamma_i-k}\right) \;\;\mbox{near $0$} \quad \mbox{and} \quad \nabla^kV_i(z) = O\left(z^{-k}\right) \;\;\mbox{near $\infty$}, \quad \forall\; k \geq 1.
\end{align}
By \eqref{e204} and \eqref{add4}, it is obvious that
$$W_k^k(z)= \frac{C_k + o(1)}{z^2}\;\; \mbox{near } 0\quad \mbox{and}\quad W_k^k(z)= O\left(z^{-2}\right)\;\; \mbox{near } \infty.$$
where $C_k$ are real constants depending on $\gamma_j$ only. Thus considering $z^2W_k^k$, we get
\BE\label{e510}
W_k^k(z)=\frac{C_k}{z^2} \quad \mbox{ in } \C.
\EE
In particular, $Z_n$ is determined uniquely. To determine $Z_k$ for $k < n$, we can do the induction step on $k$. By using \eqref{e507}, the definition of $W_k^j$, \eqref{e206} and \eqref{add4}, we obtain 
\begin{align*}
Z_k = \frac{w_k + o(1)}{z^{n+2-k}}\;\;\mbox{ near }\; 0\quad \mbox{and}\quad Z_k = O\left(\frac{1}{z^{n+2-k}}\right)\;\; \mbox{at }\;\infty,
\end{align*}
where $w_k$ is a real constant and depends only on $\gamma _j$. By the Liouville theorem, \eqref{e509} is proved.
\end{proof}

\medskip
{\it Proof of \eqref{equaf} completed}. To prove that $f$ satisfies the ODE, we use \eqref{e509} with $k=1$. By the above Lemma, for $k = 1$,
\begin{align*}
\frac{w_1}{z^{n+1}} = Z_1 = W_1^n + U_{1, z}Z_2 +  \sum_{j = 1}^{n-2}W_1^j Z_{j+2} = W_1^n + \frac{w_2}{z^n}U_{1, z} +  \sum_{j = 1}^{n-2}\frac{w_{j+2}}{z^{n - j}} W_1^j .
\end{align*}
As $f = e^{-U_1}$, we have $-U_{1, z}f = f'$ and $W_1^jf = -f^{(j+1)}$ by defintion for all $1\leq j \leq n$. Multiplying the above equation with $f$, we get 
\begin{align*}
\frac{w_1}{z^{n+1}}f = -f^{(n+1)} - \frac{w_2}{z^n}f' -  \sum_{j = 1}^{n-2}\frac{w_{j+2}}{z^{n - j}} f^{(j+1)},
\end{align*}
or equivalently
\begin{align*}
f^{(n+1)} + \sum_{k = 0}^{n-1}Z_{k+1}f^{(k)} = f^{(n+1)} + \sum_{k = 0}^{n-1}\frac{w_{k+1}}{z^{n + 1 - k}} f^{(k)} = 0.
\end{align*}
Up to change the definition of $w_k$, we are done.\qed

\subsection{Step 2}We prove that the fundamental solutions for \eqref{equaf} are just given by $f_i(z) = z^{\beta_i}$ with 
\begin{align}\label{e511}\beta_0 = -\alpha_1, \quad \beta_i = \alpha_i  - \alpha_{i+1} + i\;\;\mbox{for } 1\leq i \leq (n - 1), \quad \beta_n = \alpha_n + n.\end{align}
or equivalently we have $P(\beta_i) = 0$ where
\[
P(\beta) = \beta(\beta - 1)\ldots(\beta - n) + \sum_{i = 0}^{n - 1} w_k\beta(\beta - 1)\ldots(\beta -k + 1).
\]
By \eqref{e511}, $\beta_i$ satisfies
\BE\label{e512}
\beta_i - \beta_{i-1} = \gamma_i + 1  > 0 \quad \mbox{for all }\;1 \leq i \leq n.\EE
Let
\BE\label{e513}
f=\lambda _0|z|^{-2\alpha _1}+\displaystyle \sum _{i=1}^n \lambda_i|P_i(z)|^2,
\EE
with
\[P_i(z)=z^{(\mu _1+\mu _2+\cdots +\mu_i-\alpha _1)}+\displaystyle \sum _{j=0}^{i-1}c_{ij}z^{\mu_1+\cdots +\mu_j-\alpha _1},\]
where $\mu_i=1+\gamma _i>0$. Note that 
\[\frac{|P_i(z)|}{|z|^{\mu_1+\cdots+\mu_i-\alpha _1}}= \left|1 +\sum_{j=0}^{i-1}c_{ij}z^{-\mu_{j+1}-\cdots -\mu_i}\right| \quad \mbox{in }\; \C^*.\]
Since $|P_i(z)|$ is a single-valued function, we have $c_{ij}=0$ for $\mu_{j+1}+\cdots+\mu_i \notin \N$. In the following, we let $f^{(p,q)}$ denote $\partial _{\bar z}^q\partial_z^pf$. For any $f$ of \eqref{e513}, we define, if possible, $U=(U_1,\cdots,U_n)$ by
\begin{align}\label{e514}
e^{-U_1}=f\quad \mbox{and}\quad e^{-U_k}=2^{k(k-1)} {\rm det}_k(f) \;\; \mbox{for $2\leq k\leq n$},
\end{align}
where
\begin{align}\label{e516}
{\rm det}_k(f)={\rm det }\left(f^{(p,q)}\right)_{0\leq p,q\leq k-1} \;\; \mbox{for $1\leq k\leq n+1$}.
\end{align}

\begin{theorem}
\label{lem503}
Let ${\rm det}_k(f)$ be defined by \eqref{e516} with $f$ given by \eqref{e513} and $\lambda_i > 0$ for all $0\leq i \leq n$. Then we have ${\rm det}_k(f)>0$ in $\C^*$, $\forall\; 1\leq k\leq n$. Furthermore, $U=(U_1,...,U_n)$ defined by \eqref{e514} satisfies \eqref{e502} if and only if \eqref{equal} holds.
\end{theorem}
Before going into the details of proof of Theorem \ref{lem503}, we first explain how to construct solutions of Toda system from $f$ via the formulas \eqref{e514}. Here we follow the procedure from \cite{Y}. For any function $f$, we define ${\rm det}_k(f)$ by \eqref{e516}. Then we have
\begin{align}\label{e517}
{\rm det}_{k+1}(f)=\frac{{\rm det}_k(f)\partial _{z\bar z}{\rm det}_k(f)-\partial _z {\rm det}_k(f)\partial _{\bar z}{\rm det}_k(f)}{{\rm det}_{k-1}(f)}\quad \mbox{ for } k\geq 1.
\end{align}
The above formula comes from a general formula for the determinant of a $(k+1)\times(k+1)$ matrix. We explain it in the followings. Let ${\mathcal N} =(c_{i,j})$ be a $(k+1)\times (k+1)$ matrix:
\begin{align*}
{\mathcal N} = 
\begin{pmatrix}
M_1 &\overrightarrow{\bf u} &\overrightarrow{\bf v}\\
\overrightarrow{\bf s}& c_{k,k}&c_{k,k+1}\\
\overrightarrow{\bf t}& c_{k+1,k}&c_{k+1,k+1}
\end{pmatrix}
\end{align*}
where $\overrightarrow{\bf u}$ and $\overrightarrow{\bf v}$ stands for the column vectors consisting of first $(k-1)$ entries of the $k$-th column and $(k+1)$-th column respectively, and  $\overrightarrow{\bf s}$ and $\overrightarrow{\bf t}$ stand for rows vectors consisting of the first $(k-1)$ entries of the $k$-th rows and $(k+1)$-th rows respectively. We let
\begin{align*}
{\mathcal N}_1= \begin{pmatrix}
M_1 & \overrightarrow{\bf u}\\
\overrightarrow{\bf s} &c_{k,k}
\end{pmatrix},
\quad 
{\mathcal N}_2= \begin{pmatrix}
M_1 & \overrightarrow{\bf v}\\
\overrightarrow{\bf t} &c_{k+1,k+1}
\end{pmatrix}
\end{align*}
\begin{align*}
{\mathcal N}_1^*= \begin{pmatrix}
M_1 & \overrightarrow{\bf u}\\
\overrightarrow{\bf t} &c_{k+1,k}
\end{pmatrix},
\quad 
{\mathcal N}_2^*= \begin{pmatrix}
M_1 & \overrightarrow{\bf v}\\
\overrightarrow{\bf s} &c_{k,k+1}
       \end{pmatrix}.
\end{align*}
Then we have
\[{\rm det}({\mathcal N}){\rm det}(M_1) = {\rm det}({\mathcal N}_1){\rm det}({\mathcal N}_2) - {\rm det}({\mathcal N}_1^*){\rm det}({\mathcal N}_2^*).\]
Since the proof is elementary, we omit it. Clearly, \eqref{e517} follows from the above formula immediately.

\medskip
Suppose that ${\rm det}_k(f) > 0$ for $1\leq k \leq n$ and ${\rm det}_{n+1}(f)=2^{-n(n+1)}$. Define $U_1$ by $f=e^{-U_1}$. As $-e^{-2U_1}U_{1,z\bar z}=ff_{z\bar z}-f_zf_{\bar z}$, then
\[-4U_{1,z\bar z}=e^{2U_1-U_2} \quad \mbox{if and only if} \quad e^{-U_2} = 4(ff_{z\bar z}-f_zf_{\bar z})=4{\rm det}_2(f).\]
By the induction on $k$, $2\leq k\leq n$, we have 
\begin{align*}
-4e^{-2U_k}U_{k,z\bar z}&=4e^{-2U_k}\big[\log {\rm det}_k(f)\big]_{z\bar z}\\
&=4\cdot 2^{2k(k-1)}\big[{\rm det}_k(f)\partial _{z\bar z}{\rm det}_k(f)-\partial _z {\rm det}_k(f)\partial _{\bar z}{\rm det}_k(f)\big]\\
&=2^{2k(k-1)+2}{\rm det}_{k+1}(f)\;{\rm det}_{k-1}(f)\\
&=2^{(k+1)k}e^{-U_{k-1}}{\rm det}_{k+1}(f).
\end{align*}
Thus, $U_k$ satisfies $\triangle U_{k,z\bar z}+e^{2U_k - U_{k+1}-U_{k-1}}=0$ in $\C^*$ if and only if $e^{-U_{k+1}}=2^{(k+1)k}{\rm det}_{k+1}(f)$. For the last equation $k=n$, we have
\[-4e^{-2U_n}U_{n,z\bar z} = 2^{(n+1)n}e^{-U_{n-1}}{\rm det}_{n+1}(f).\]
Thus, $U_n$ satisfies $\triangle U_n+e^{2U_n-U_{n-1}}=0$ in $\C^*$ if and only if ${\rm det}_{n+1}(f) =2^{-n(n+1)}$.

\bigskip
Therefore, assume that $U = (U_k)$ given by \eqref{e514}, \eqref{e516} and \eqref{e513} is a solution of the Toda system \eqref{e502}, to get the equality in \eqref{equal}, it is equivalently to show
\begin{align}
\label{equafs}
{\rm det}_{n+1}(f) = \lambda_0\lambda_1\cdots\lambda_n \times \Pi_{1 \leq i\leq j\leq n}\left(\sum_{k= i}^j\mu_k\right)^2
\end{align}
for $f$ given by \eqref{e513}. We have first
\begin{lemma}\label{lem504}
Let $g = |z|^{2\beta}f$ with $\beta \in \R$, and $f$ be a complex analytic function in $\C^*$, there holds
\begin{align}
\label{zbeta}
{\rm det}_k(g) = |z|^{2k\beta}{\rm det}_k(f)  \;\;\mbox{in } \; \C^*, \quad \forall\; k \in \N^*.
\end{align}
\end{lemma}
\begin{proof} This is obviously true for $k = 1$, we can check also easily for $k = 2$. Suppose that the above formula holds for $1 \leq \ell \leq k$, then by formula \eqref{e517}, 
\begin{align*}
{\rm det}_{k+1}(g) = \frac{{\rm det}_k(g)\p_{z\bar z}{\rm det}_k(g) - \p_{z}{\rm det}_k(g)\p_{\bar z}{\rm det}_k(g)}{{\rm det}_{k-1}(g)} & = \frac{{\rm det}_2\left({\rm det}_k(g)\right)}{{\rm det}_{k-1}(g)}\\
& = \frac{{\rm det}_2\left(|z|^{2k\beta}{\rm det}_k(f) \right)}{|z|^{2(k-1)\beta}{\rm det}_{k-1}(f) }\\
& = |z|^{2(k+1)\beta}\frac{{\rm det}_2\left({\rm det}_k(f) \right)}{{\rm det}_{k-1}(f)}\\
& = |z|^{2(k+1)\beta}{\rm det}_{k+1}(f).
\end{align*}
The equality \eqref{zbeta} holds when ${\rm det}_{k-1}(f) \ne 0$.\end{proof}

\medskip
Thanks to \eqref{zbeta}, to prove \eqref{equafs}, it is enough to prove the following: Let 
\begin{align}
\label{tildef}
\widetilde f = \lambda_0 + \sum_{i=1}^n \lambda_i |P_i(z)|^2  \quad\mbox{in }\C
\end{align}
with $P_i$ given by \eqref{pi}, then
\begin{align}
\label{det}
{\rm det}_{n+1}(\widetilde f) = \lambda_0\lambda_1\cdots\lambda_n \times \prod_{1 \leq i\leq j\leq n}\left(\sum_{k= i}^j\mu_k\right)^2\times |z|^{2n\gamma_1 + 2(n-1)\gamma_2 +  \ldots + 2 \gamma_n}.
\end{align}
Here we used $(n+1)\alpha_1 = n\gamma_1 + (n-1)\gamma_2 +  \ldots + \gamma_n$ for $SU(n+1)$. 

\medskip
{\it Proof of \eqref{det}}. We proceed by induction. Let $n = 1$, we have $P_1 = c_0 + z^{\mu_1}$, so
$${\rm det}_2(\widetilde f) = {\rm det}_2\left(\lambda_0 + \lambda_1|P_1|^2\right) = |z|^{-4\alpha_1}\lambda_0\lambda_1|P_1'|^2 = \lambda_0\lambda_1\mu_1^2|z|^{2(\mu_1 - 1)} = \lambda_0\lambda_1\mu_1^2 |z|^{2\gamma_1}.$$
since $\mu_1 - 1 = \gamma_1$. Then \eqref{det} holds true for $n = 1$. 

\medskip
Suppose that \eqref{det} is true for some $(n-1) \in \N^*$, we will prove \eqref{det} for the range $n$. Define $L_k(P)$ to be the vertical vector $(P, \p_z P, \ldots, \p^k_zP) \in \C^{k+1}$ for any smooth function $P$ and $k \in \N^*$. Denote $P_0 \equiv 1$, there holds
\begin{align*}
{\rm det}_{n+1}(\widetilde f) & = \sum_{0 \leq i_k\leq n, i_p \ne i_q} \lambda_{i_0}\lambda_{i_1}\cdots\lambda_{i_n}{\rm det}\Big(\overline{P_{i_0}}L_n(P_{i_0}), \p_{\bar z}\overline{P_{i_1}}L_n(P_{i_1}), \cdots, \p^n_{\bar z}\overline{P_{i_n}}L_n(P_{i_n})\Big)\\
& = \lambda_0\lambda_1\cdots\lambda_n\sum_{1 \leq i_k\leq n, i_p \ne i_q} {\rm det}\Big(\overline{P_0}L_n(P_0), \p_{\bar z}\overline{P_{i_1}}L_n(P_{i_1}), \cdots, \p^n_{\bar z}\overline{P_{i_n}}L_n(P_{i_n})\Big).
\end{align*}
The last line is due to $P_0 \equiv 1$. Let $e_1$ be the vertical vector $(1, 0, \ldots, 0)$, we have 
\begin{align*}
& {\rm det}\Big(\overline{P_0}L_n(P_0), \p_{\bar z}\overline{P_{i_1}}L_n(P_{i_1}), \cdots, \p^n_{\bar z}\overline{P_{i_n}}L_n(P_{i_n})\Big)\\
= & \;{\rm det}\Big(e_1, \p_{\bar z}\overline{P_{i_1}}L_n(P_{i_1}), \cdots, \p^n_{\bar z}\overline{P_{i_n}}L_n(P_{i_n})\Big)\\
= & \;{\rm det}\Big(\overline{P_{i_1}'}L_{n - 1}(P_{i_1}'), \cdots,\overline{P_{i_n}'}L_{n - 1}(P_{i_n}')\Big).
\end{align*}
Therefore ${\rm det}_{n+1}(\widetilde f) = \lambda_0\lambda_1\cdots\lambda_n{\rm det}_n(h)$
with $h = \sum_{1\leq i \leq n}|P_i'|^2$. Moreover, for $i \geq 1$,
\begin{align*}P_i' & = \sum_{k = 1}^{i - 1}(\mu_1+\mu_2+\ldots+\mu_k)c_{ik}z^{\mu_1+\mu_2+\ldots+\mu_k - 1} + (\mu_1+\mu_2+\ldots+\mu_i)z^{\mu_1+\mu_2+\ldots+\mu_i - 1}\\
& = (\mu_1+\mu_2+\ldots+\mu_i)z^{\mu_1 - 1}\widetilde P_i
\end{align*}
where
$$\widetilde P_i = z^{\mu_2+\ldots+\mu_i} + \sum_{k = 1}^{i - 1}\widetilde c_{ik}z^{\mu_2+\ldots+\mu_k}\quad \mbox{with }\; \widetilde c_{ij} \in \C.$$
This means that
\begin{align*}h & = |z|^{2\gamma_1}\Big[\sum_{i = 1}^n (\mu_1+\mu_2+\ldots+\mu_i)^2|\widetilde P_i|^2\Big]\\
& = |z|^{2\gamma_1}\Big[ \mu_1^2 +  \sum_{i = 1}^{n - 1}(\mu_1+\mu_2+\ldots+\mu_{i+1})^2|\widetilde P_{i+1}|^2\Big] := |z|^{2\gamma_1}\widetilde h,\end{align*}
hence $\widetilde h$ is in the form of \eqref{tildef} with $(n - 1)$. Consequently, by the induction hypothesis, we get 
\begin{align*}
{\rm det}_{n+1}(\widetilde f) = & \;\lambda_0\lambda_1\cdots\lambda_n {\rm det}_n(h)\\ = & \; \lambda_0\lambda_1\cdots\lambda_n |z|^{2n\gamma_1}{\rm det}_n(\widetilde h)\\
= & \;\lambda_0\lambda_1\cdots\lambda_n |z|^{2n\gamma_1}\\
& \; \times \prod_{1 \leq k \leq n}(\mu_1 + \mu_2 + \ldots \mu_k)^2\times \prod_{2 \leq i\leq j\leq n}\left(\sum_{k= i}^j\mu_k\right)^2\times |z|^{2(n-1)\gamma_2 +  \ldots + 2 \gamma_n},
\end{align*}
which yields clearly the equality \eqref{det}. \qed

\medskip
On the other hand, assume that \eqref{equal} holds true, using the above analysis and \eqref{equafs}, we see that $U$ defined by \eqref{e514} and \eqref{e513} is a solution of \eqref{e502} in $\C^*$ provided that ${\rm det}_k(f)>0$ in $\C^*$. 

\medskip
First we make a general calculus of ${\rm det}_k(g)$ with 
\begin{align}
g = \sum_{i, j = 0}^n m_{ij}f_i\overline{f_j}, \quad \mbox{where $m_{ij} = \overline{m_{ji}}$ for all $0 \leq i, j \leq n$,}
\end{align}
where $f_i(z)=z^{\beta _i}$. Let $M = (m_{ij})_{0\leq i, j\leq n}$ and $J = (z_{ij})_{0\leq i, j\leq n}$ with $z_{ij} = \left(z^{\beta_j}\right)^{(i)}$. Let ${\mathcal N}_{i_1, \ldots, i_k}^{j_1, \ldots, j_k}$ be the $k\times k$ sub matrix $(b_{ij})_{i = i_1, \ldots, i_k, j= j_1, \ldots, j_k}$, for any matrix ${\mathcal N} = (b_{ij})$, we denote also  ${\mathcal N}_{i_1, \ldots, i_k}$ the $k \times (n+1)$ sub matrix by taking the rows $i_1, \ldots, i_k$ of ${\mathcal N}$, and ${\mathcal N}^t$ means the transposed matrix of ${\mathcal N}$. 

\medskip
As $g^{(p, q)} = \sum m_{ij}f_i^{(p)}\overline{f_j^{(q)}}$. For $1 \leq k \leq n$, we can check easily that
\begin{align*}
\left(g^{(p, q)}\right)_{0\leq p,q \leq k} = J_{0, 1, \ldots, k}M \overline{J_{0, 1, \ldots, k}}^t,
\end{align*}
and
\begin{align}
\label{add5}
\begin{split}
& \;  {\rm det}\left(J_{0, 1, \ldots, k}M \overline{J_{0, 1, \ldots, k}}^t\right)\\
= & \; \sum_{0\leq i_0 < i_1 < \ldots<  i_k \leq n, 0 \leq j_0 < j_1 \ldots < j_k \leq n}{\rm det}\left(J_{0, 1, \ldots, k}^{i_0, i_1, \ldots i_k}M_{i_0, i_1, \ldots i_k}^{j_0, j_1\ldots, j_k}\overline{J_{0, 1, \ldots, k}^{j_0, j_1\ldots, j_k}}^t\right)\\
= & \; \sum_{0\leq i_0 < i_1 < \ldots<  i_k \leq n, 0 \leq j_0 < j_1 \ldots < j_k \leq n}{\rm det}\left(M_{i_0, i_1, \ldots i_k}^{j_0, j_1\ldots, j_k}\right){\rm det}\left(J_{0, 1, \ldots, k}^{i_0, i_1, \ldots i_k}\right)\overline{{\rm det}\left(J_{0, 1, \ldots, k}^{j_0, j_1\ldots, j_k}\right)}.
\end{split}
\end{align}
Moreover, exactly as for \eqref{det}, by induction, we can prove that
\begin{align}
\label{add6}
\begin{split}
{\rm det}\left(J_{0, 1, \ldots, k}^{i_0, i_1, \ldots i_k}\right) & = \prod_{0 \leq p < q\leq k}\left(\beta_{i_q} - \beta_{i_p}\right)\times z^{(k+1)\beta_{i_0} + k(\beta_{i_1} - \beta_{i_0} - 1) + \ldots + (\beta_{i_k} - \beta_{i_{k-1}} - 1)}\\
& =  \prod_{0 \leq p < q\leq k}\left(\beta_{i_q} - \beta_{i_p}\right)\times z^{\beta_{i_0} + \beta_{i_1} + \ldots + \beta_{i_k} - \frac{k(k+1)}{2}}.
\end{split}
\end{align}

Given $f$ by \eqref{e513} with $\lambda_i$ satisfying \eqref{equal}, we will prove that ${\rm det}_k(f)>0$ in $\C^*$. Clearly, $f > 0$ in $\C^*$ and $f = \sum_{0\leq i, j \leq n} m_{ij}f_i\overline{f_j}$ where
\begin{align*}
M = (m_{ij}) = B \overline{B}^t, \quad B = (b_{ij}) \;\; \mbox{with }\; b_{ii} = \sqrt{\lambda_i}, \; b_{ij} = \sqrt{\lambda_i}c_{ji} \;\mbox{for $j > i$}, \; b_{ij} = 0 \;\mbox{for $j < i$}.
\end{align*}
For $1 \leq k \leq n$, denote ${\mathcal B} = J_{0, 1, \ldots, k}B$, we can check that
\begin{align*}
{\rm det}_{k+1}(f) = {\rm det}\left(J_{0, 1, \ldots, k}M \overline{J_{0, 1, \ldots, k}}^t\right) & = {\rm det}\left({\mathcal B}\overline{{\mathcal B}}^t\right)\\
& =  \sum_{0\leq i_0 < i_1 < \ldots<  i_k \leq n} {\rm det}\left({\mathcal B}_{0, 1, \ldots, k}^{i_0, i_1, \ldots i_k}\right){\rm det}\left(\overline{{\mathcal B}_{0, 1, \ldots, k}^{i_0, i_1\ldots, i_k}}^t\right)\\
& =  \sum_{0\leq i_0 < i_1 < \ldots<  i_k \leq n} \left|{\rm det}\left({\mathcal B}_{0, 1, \ldots, k}^{i_0, i_1, \ldots i_k}\right)\right|^2.
\end{align*}
As ${\rm det}_{n+1}(f) = 2^{-n(n+1)} \ne 0$ by \eqref{equal} and \eqref{equafs}, the rank of the matrix ${\mathcal B}$ must be $(k + 1)$ in $\C^*$, hence for any $z \in \C^*$, we have $0\leq i_0 < i_1 < \ldots<  i_k \leq n$, such that ${\rm det}\left({\mathcal B}_{0, 1, \ldots, k}^{i_0, i_1, \ldots i_k}\right)(z) \ne 0$, thus ${\rm det}_{k+1}(f)  > 0$ in $\C^*$.

\medskip
To complete the proof of Theorem \ref{lem503}, it remains to compute the strength of the singularity. Notice that $M = B \overline{B}^t$ is a positive hermitian matrix, since $\lambda_i > 0$. By the formulas \eqref{add5}, \eqref{add6}, as $i_p \geq p$, $j_p \geq p$, $\beta_i$ are increasing and 
$$\sum_{p = 0}^k \beta_p - \frac{k(k+1)}{2} = -(k+1)\alpha_0 + k\gamma_1 + (k-1)\gamma_2 +\ldots + \gamma_k = -\alpha_{k+1},$$
we get
\begin{align}
\label{add7}
{\rm det}_{k+1}(f) = \prod_{0\leq p<q\leq k}(\beta _q-\beta _p) |z|^{-2\alpha _{k+1}}\big[\zeta_k + o(1)\big] \quad \mbox{as $z\to 0$},
\end{align}
with $\zeta_k = {\rm det}\left(M_{0, 1, \ldots, k}^{0, 1, \ldots, k}\right) > 0$. This implies 
$$U_{k+1}= -2\alpha _{k+1}\log |z|+O(1) \quad \mbox{ near $0$}.$$ Hence $U = (U_1,\cdots,U_n)$ satisfies \eqref{e502} in $\C$. This completes the proof of Theorem \ref{lem503}. \qed

\medskip
By Theorem \ref{lem503}, we have proved that any $f$ given by \eqref{e513} verifying \eqref{equal} is a solution of \eqref{equaf}, because $U=(U_1,\ldots,U_n)$ defined by \eqref{e514} is a solution of the Toda system. In particular, it is the case for $f= \sum _{0\leq i\leq n}\lambda _i|z|^{2\beta _i}$ satisfying \eqref{equal}, with $\beta_i$ are given by \eqref{e511}. Let $L$ denote the linear operator of the differential equation \eqref{equaf}. Then
\[0=\overline {L}L(f)=\displaystyle \sum_{i=0}^n\lambda _i|L(z^{\beta _i})|^2,\]
which implies $L(z^{\beta _i})=0$, $\forall\; 0\leq i\leq n$. Thus Step 2 is proved.

\subsection{Step 3}
Suppose $U=(U_1,\ldots, U_n)$ is a solution of equation \eqref{e502}, we will prove that $f=e^{-U_1}$ can be  written as the form of \eqref{e513}. For any solution $(U_i)$, as $f = e^{-U_1} > 0$ satisfies \eqref{equaf}, we have
\begin{align*}
f = \sum_{i, j = 0}^n m_{ij}f_i\overline{f_j}, \quad \mbox{where $m_{ij} = \overline{m_{ji}}$ for all $0 \leq i, j \leq n$,}
\end{align*}
where $f_i(z)=z^{\beta _i}$ is a set of fundamental solutions of \eqref{equaf}.

\medskip
We want to prove that $f$ can be written as a sum of $|P_i(z)|^2$, which is not true in general, because even a positive polynomial in $\C$ cannot be written always as sum of squares of module of polynomials. For example, it is the case for $2|z|^6 - |z|^4 - |z|^2 + 2$. It means that, we need to use further informations from the Toda system. In fact, we will prove that $M = (m_{ij})$ is a positive hermitian matrix.

\medskip
With $V_i$ given by \eqref{add3},
$$e^{V_1} = |z|^{2\alpha _1} e^{-U_1} = |z|^{2\alpha_1}f = m_{00} + \sum_{i = 1}^n m_{ii} |z|^{2(\beta_i - \beta_0)} + 2\sum_{0 \leq i < j \leq n}{\rm Re}\left(m_{ij}{\bar z}^{\beta_j- \beta_i}\right)|z|^{2(\beta_i-\beta_0)},$$
Take $z = 0$, we get $m_{00} > 0$. Let $J = (z_{ij})_{0\leq i, j\leq n}$ with $z_{ij} = \left(z^{\beta_j}\right)^{(i)}$ as in Step 2. Using \eqref{add5}, \eqref{add6} and the monotonicity of $\beta_i$, exactly as before, we get, for $1 \leq k \leq n - 1$ 
\begin{align*}
{\rm det}_{k+1}(f) = \prod_{0\leq p<q\leq k}(\beta _q-\beta _p) |z|^{-2\alpha _{k+1}}\left[{\rm det}\left(M_{0, 1, \ldots, k}^{0, 1, \ldots, k}\right) + o(1)\right], \quad \mbox{as $z\to 0$}.
\end{align*}
Recall that $e^{-U_{k+1}} = 2^{k(k+1)}{\rm det}_{k+1}(f)$ and $V_{k+1}$ is defined by \eqref{add3}, 
\begin{align*}
\frac{e^{-V_{k+1}(0)}}{2^{2(k+1)k}} = \left[|z|^{2\alpha_{k+1}}{\rm det}_{k+1}(f)\right]_{z = 0} = {\rm det}\left(M_{0,1, \ldots k}^{0, 1\ldots, k}\right)\times \prod_{0 \leq p < q\leq k}\left(\beta_q - \beta_p\right)^2,
\end{align*}
which yields
\begin{align}
\label{signMk}
{\rm det}\left(M_{0,1, \ldots k}^{0, 1\ldots, k}\right) > 0, \quad \forall\; 1 \leq k \leq n-1.
\end{align}

Similarly, when $k = n$, noticing that 
$$\sum_{p = 0}^n \beta_p - \frac{n(n+1)}{2} = 0,$$
we obtain 
\begin{align}
\label{detMn}
2^{-n(n+1)} = {\rm det}_{n+1}(f) = {\rm det}(M)\times \prod_{0 \leq p < q\leq n}\left(\beta_q - \beta_p\right)^2,
\end{align}
hence ${\rm det}(M) > 0$. Combining with \eqref{signMk} and $m_{00} > 0$, it is well known that $M$ is a positive hermitian matrix. Consequently, we can decompose $M = B\overline{B}^t$ with a upper triangle matrix $B = (b_{ij})$ where $b_{ii} > 0$. To conclude, we have 
\begin{align*}
f = \sum_{i, j = 0}^n m_{ij}f_i\overline{f_j} = \sum_{k= 0}^n |Q_k|^2, \quad \mbox{where }\; Q_k = \sum_{i= 0}^k b_{ik}f_i.
\end{align*}
It is equivalent to saying that $f$ is in the form of \eqref{e513} with $\lambda_i = b_{ii}^2 > 0$. Combining with Theorem \ref{lem503}, the proof of Theorem \ref{mainsingular} is finished.\qed

\section{Quantization and Nondegeneracy}
\reset
Here we will prove Theorem \ref{thm105}. We first prove the quantization of the integral of $e^{u_i}$. By \eqref{add5}, \eqref{add6} and again the monotonicity of $\beta_i$ with $f$ given by \eqref{e513}, we have for $1\leq k \leq n$,
\[e^{-U_k}=2^{k(k-1)}{\rm det}_k(f)=|z|^{2(\beta_{n-k+1}+\cdots+\beta _n)-k(k-1)}\left[c_k+ o(1)\right], \quad \mbox{as $|z| \to \infty$},\]
where 
$$c_k = 2^{k(k-1)}\lambda_{n-k+1}\lambda_{n-k+2}\cdots\lambda_n \times \prod_{n-k+1\leq q < p \leq n}(\beta_p - \beta_q)^2 >0.$$ 
Thus, as $-\Delta U_k = e^{u_k} - 4\pi\alpha_k\delta_0$
\begin{align*}
\int_{\R^2}e^{u_k}dx &= 4\pi \alpha _k+\lim _{R\rightarrow +\infty}\int _{\partial B_R}\frac{\partial U_k}{\partial \nu}ds\\
&= 4\pi\left[\alpha_k+\beta_{n-k+1} + \cdots +\beta_n-\frac{k(k-1)}{2}\right]\\
&= 4\pi\Big[\alpha_k+\alpha_{n-k+1}+k(n-k+1)\Big].
\end{align*}
Therefore,
\[\sum _{k=1}^n a_{ik}\int_{\R^2}e^{u_k}dx=4\pi(2+\gamma_k+\gamma_{n+1 - k}),\]
which implies
\[u_k(z)=-4\pi (2 + \gamma_{n+1 - k})\log |z|+O(1), \quad \mbox{for large $|z|$.}\]
This proves the quantization.

\medskip
To prove the nondegeneracy, we let $(u_i)$ be a solution of the singular Toda system $SU(n+1)$ \eqref{e103n} and $\phi_i$ be solutions of the linearized system $LSU(n+1)$:
\BE\label{b1} - \Delta \phi_i = \D \sum_{j = 1}^n a_{ij}e^{u_j}\phi_j \;\; \text{in }\mathbb R^2, \quad
\phi_i \in L^\infty(\R^2)\;\; \forall\; 1\leq i \leq n, \EE 
or equivalently 
\begin{align*}
-4\Phi_{i, z\bar z} = {\rm exp}\left({\sum_{j = 1}^n
a_{ij}U_j}\right)\times \sum_{j = 1}^n \left(a_{ij}\Phi_j\right)\;\;
\text{in }\R^2, \quad
\Phi_i \in L^\infty(\R^2), \quad \forall\; 1\leq i \leq n
\end{align*}
where $U_j$ are defined by \eqref{singSU} and $\Phi_j$ defined by \eqref{b2}. 

\medskip
We will use the quantities $Y_1^j =
e^{U_1}\left[\left(e^{-U_1}\Phi_1\right)^{(j+1)} -
\left(e^{-U_1}\right)^{(j+1)}\Phi_1\right]$ for $1 \leq j \leq n$, and  
\begin{align*} 
Y_{k+1}^j = - \frac{Y_{k, \bar z}^j + W_{k+1}^j\Phi_{k, z\bar
z}}{U_{k, z\bar z}}\quad \mbox{for $1\leq k < j \leq n$}.\end{align*}
Recall that $Y_{n, \bar z}^n = 0$ in $\C^*$ for solutions of $LSU(n+1)$, we can prove also (as for \eqref{e205})
\begin{align}
\label{Yjj}
Y_{j, \bar z}^j = -\Phi_{j, z\bar z}U_{j+1, z} - U_{j, z\bar z}\Phi_{j+1, z} \quad \mbox{for solutions of $LSU(n+1)$ and $j < n$}.
\end{align}

Now we define some new invariants $\widetilde Z_k$ for solutions of \eqref{b1}, which correspond to $Z_k$ for system $SU(n+1)$. Let
\begin{align*}
\widetilde Z_n = Y_n^n, \quad \mbox{and}\quad \widetilde Z_k = Y_k^n + \Phi_{k, z}Z_{k +1} + \sum_{j = k}^{n-2}Y_k^j Z_{j+2}, \;\; \forall\;  k = n-1, n-2, \ldots, 1.
\end{align*}
The central argument is 
\begin{lemma}
For any solution of \eqref{b1}, we have $\widetilde Z_k \equiv 0$ in $\C^*$ for all $1 \leq k \leq n$.
\end{lemma}

\proof By the same argument as in section 4, we have that $\widetilde Z_n$ is holomorphic in $\C^*$, since 
\begin{align*}
\widetilde Z_n = Y_n^n = \sum_{i = 1}^n \Phi_{i, zz} - 2\sum_{i = 1}^n
U_{i, z}\Phi_{i, z} + \sum_{i = 1}^{n-1} \left(\Phi_{i, z}U_{i+1, z}
+ U_{i, z}\Phi_{i+1, z}\right).
\end{align*}
Using the integral representation formula for $\Phi_i$, we see that $\nabla^k\Phi_i = O(z^{-k})$ as $|z|\to \infty$ for all $k \geq 1$, so $\widetilde Z_n = O(z^{-2})$ at infinity. On the other hand, since $\gamma_j > -1$ for all $1\leq j \leq n$, we have $\Phi_i \in C^{0, \alpha}(\C)$ with some $\alpha \in (0, 1)$, for any $1\leq i \leq n$. Again, by elliptic estimates, we can claim that 
\begin{align*}
\nabla^k\Phi_i(z) = o\left(z^{-k}\right) \;\; \mbox{as $z \to 0$}, \quad \mbox{for $k \geq 1$, $1\leq i\leq n$.}
\end{align*}
By the behavior of $U_i$ via \eqref{add4}, $\widetilde Z_n = o(z^{-2})$ near the origin, so $\widetilde Z_n \equiv 0$ in $\C^*$. 

\medskip
Combining the iterative relations on $Y_k^j$, the behaviors of $\Phi_i$ and $U_j$, we can claim that for all $k \leq j \leq n$,
\begin{align}
\label{asymY}
Y_k^j = O\left(z^{k-j-2}\right) \;\;\mbox{as $|z|\to \infty$}\quad \mbox{and} \quad Y_k^j = o\left(z^{k-j-2}\right) \;\;\mbox{as $|z|\to 0$}.
\end{align}
Therefore (recalling that $Z_k = w_kz^{k-2-n}$ for any $k$), as $Z_n = W_n^n$ and $Y_n^n = 0$,
\begin{align*}
\widetilde Z_{n-1, \bar z} = Y_{n-1, \bar z}^n + \Phi_{n-1, z\bar z}Z_n = -U_{n-1, \bar z}Y_n^n - \Phi_{n-1, z\bar z}W_n^n + \Phi_{n-1, z\bar z}Z_n = 0.
\end{align*}
So $\widetilde Z_{n-1}$ is holomorphic in $\C^*$. Using expression of $Z_k$, the asymptotic behavior of $\Phi_i$ and \eqref{asymY}, we see that $\widetilde Z_{n-1} = O(z^{-3})$ at infinity and $\widetilde Z_{n-1} = o(z^{-3})$ near $0$, hence $\widetilde Z_{n-1} = 0$ in $\C^*$. For $k \leq n - 2$, suppose that $\widetilde Z_j = 0$ for $j > k$, we have
\begin{align*}
\widetilde Z_{k, \bar z} = & \; Y_{k, \bar z}^n + \Phi_{k, z\bar z}Z_{k+1} + Y_{k, \bar z}^kZ_{k+2} + \sum_{j = k+1}^{n-2} Y_{k, \bar z}^jZ_{j+2}\\
= & \; -U_{k, z\bar z}\Big[Y_{k+1}^n + \Phi_{k+1, z}Z_{k+2} + \sum_{j = k+1}^{n-2} Y_{k+1}^jZ_{j+2}\Big]\\
& + \Phi_{k, z\bar z}\Big[Z_{k+1} - W_{k+1}^n - U_{k+1, z}Z_{k+2} - \sum_{j = k+1}^{n-2} W_{k+1}^jZ_{j+2}\Big]\\
= & \; -U_{k, z\bar z}\widetilde Z_{k+1}\\
= & \; 0.
\end{align*}
Here we used the definition of $Z_{k+1}$. Similarly, the asymptotic behaviors yield that $\widetilde Z_k = 0$ in $\C^*$. The backward induction finishes the proof. \qed

\medskip
Let $g = f\Phi_1$ with $f = e^{-U_1}$, by the definition of $Y_1^j$, we see that $g^{(j+1)} = f^{(j+1)}\Phi_1 + fY_1^j$ for any $1\leq j \leq n$. Finally, 
\begin{align*}
g^{(n+1)} = f^{(n+1)}\Phi_1 + fY_1^n & = -\Phi_1\sum_{j = 0}^{n - 1}Z_{j+1}f^{(j)} + fY_1^j\\
& = -Z_1f\Phi_1 - Z_2f'\Phi_1 - \sum_{j = 2}^{n - 1}Z_{j+1}\left[g^{(j)} - fY_1^{j-1}\right] + fY_1^n\\
& = -\sum_{j = 0}^{n - 1}Z_{j+1}g^{(j)} + f\Big[Y_1^n + \Phi_{1, z}Z_2 - \sum_{j = 1}^{n - 2}Y_1^jZ_{j+2}\Big]\\
& = -\sum_{j = 0}^{n - 1}Z_{j+1}g^{(j)}.
\end{align*}
For the last line, we used $\widetilde Z_1 = 0$. Therefore $g$ satisfies exactly the same differential equation \eqref{equaf} for $f$. 

\medskip
As $g$ is a real function in $\C^*$, we get $g = \sum \widetilde m_{kl}f_k\overline{f_l}$ with a hermitian matrix $(\widetilde m_{kl})$. As before, the coefficients $\widetilde m_{kl}$ need to be zero if $\mu_{k+1}+\cdots+\mu_l \notin \N$, $k<l$, because for $z=|z|e^{i\theta}$, 
\[g=\displaystyle \sum_{k=0}^n \widetilde{m}_{kk}|z|^{2\beta _k}+ 2\sum_k|z|^{2\beta _k} {\rm Re}\left(\sum _{k<l}\widetilde{m}_{kl}e^{i(\mu_{k+1}+\cdots+\mu_l)\theta}\right) \] 
is a single-valued function in $\C^*$. Besides, we can also eliminate the subspace of constant functions for $\Phi_1$ as in section 4. We can conclude then the solution space for \eqref{b1} has the same dimension for the solution manifold for \eqref{e103n}, which means just the nondegeneracy. \qed

\section{Proof of Theorem \ref{thm101}}
\reset
Let $u$ be a solution of \eqref{a4}. By the proof of Lemma \ref{lem502n}, $f=e^{-U_1}$ satisfies the differential equation:
\begin{align}\label{e701}
L(f)=f^{(n+1)}+\sum_{k=0}^{n-1}Z_{k+1}f^{(k)} =0 \quad \mbox{in} \; \C\setminus \{P_1,\ldots,P_m\},
\end{align}
where $Z_{k+1}$ is a meromorphic function with poles at $\{P_1,\ldots,P_m\}$ and $Z_{k+1}(z)=O(|z|^{-n+k-1})$ at $\infty$. 

\noindent From Lemma \ref{lem21}, the principal part of $Z_k$ at $P_j$ is
\begin{align}\label{e702}
Z_k=\frac{w_k}{(z-P_j)^{n+1-k}}+O\left(\frac{1}{|z-P_j|^{n-k}}\right),
\end{align}
where the coefficient depends only on $\left\{ \gamma _{ij}, 1\leq i \leq n\right\}$.

\medskip
As we knew in the Introduction, locally $f$ can be written as a sum of $|\nu _i(z)|^2$, where $\nu_i(z)$ is a holomorphic function. Hence
\[0=\overline{L}L(f)=\sum _{i=0}^n|L(\nu _i)|^2\]

\noindent Therefore, $\{\nu_i\}_{0\leq i \leq n}$ is a set of fundamental solutions of \eqref{e701}, and by \eqref{e701}, $\left\|\nu \wedge \cdots \wedge \nu^{(n)}(z)\right\|$ remains a constant through its analytical continuation. The local exponents $\left\{\beta _{ij}, 1\leq i \leq n\right\}$ of \eqref{e701} at each $P_j$ is completely determined by the principal part of $Z_k$. Hence by \eqref{e702} and \eqref{e511}, we have 
\[\beta_{0j}=-\alpha_{1j},\quad \beta_{ij}=\beta_{i-1,j}+\gamma_{ij}+1.\] 
Therefore, near each $P_\ell$, $\ell=1,2,\ldots,m$, $\nu_i(P_\ell +z) =\sum _{0\leq j\leq n} c_{ij} z^{\beta _{j\ell}}g_j(z)$, where $g_j$ is a holomorphic function in a neighborhood of $P_\ell$. Since $\beta_{j\ell}-\beta_{0\ell}$ are positive integers, we have
\begin{align}\label{e703}
\nu(P_\ell+ze^{2\pi i})=e^{2\pi i\beta _{0\ell}}\nu (P_\ell+z),
\end{align}
i.e.~the monodromy of $\nu$ near $P_\ell$ is $e^{2\pi i\beta _{0\ell}}I$, $I$ is the identity matrix. Therefore, the monodromy group of \eqref{equaf} consists of scalar multiples of $I$ only, which implies $[\nu (z)]$, as a map into $\mathbb{CP}^n$, is smooth at $P_\ell$ and well-defined in $\C$.

\medskip
Applying the estimate of Brezis and Merle \cite{BM}, we have
\[u_i(z)=-(4+2\gamma _i^*)\log |z|+O(1)\quad \mbox{at } \infty,\]
for some $\gamma_i^*$. To compute $\gamma_i^*$, we might use the Kelvin transformation, $\widehat{u}_i(z)=u_i(z|z|^{-2})-4\log |z|$. Then $\widehat{u}_i(z)$ also satisfies \eqref{a4} with a new singularity at $0$,
\[\widehat{u}_i(z)=-2\gamma_i^*\log|z|+O(1)\quad \mbox{near }0.\]

The local exponent of ODE \eqref{e701} corresponding to $\widehat{u}_i$ near $0$ is $\beta_i^*$ where $\beta_i^*-\beta_{i-1}^*=\gamma_i^*+1$ for $1\leq i\leq n$. Let $\widehat{\nu} = (\widehat{\nu}_1,\cdots,\widehat{\nu}_n)$ be a holomorphic curve corresponding to $\widehat{u}$, then
\[\widehat{\nu}_i(ze^{2\pi i})=e^{2\pi i \beta _i^*}\widehat{\nu}_i(z).\] 
Since the monodromy near $0$ is a scalar multiple of the identity matrix, we conclude that $\beta_i^*-\beta_0^*$ must be integers and therefore, all $\gamma_i^*$ are integers. By identifying $S^2=\C\cup \{\infty \}$, we see $\nu(z)$ can be smoothly extended to be a holomorphic curve from $S^2$ into $\mathbb{CP}^n$ and $\infty$ might be a ramificated point with the total ramification index $\gamma_i^*$. This ends the proof of Theorem \ref{thm101}. \qed

\bigskip

\section{Appendix: explicit formula for $SU(3)$}
\reset
For general $SU(n+1)$ Toda system \eqref{e103n}, depending the valus of $\gamma_i > -1$, we can have many different situations by Theorem \ref{thm101n}. The solution manifolds have  dimensions ranging from $n$ to $n(n+2)$. On the other hand, with the expression of $U_1$ given by \eqref{e107n} and $f = e^{-U_1}$, we can obtain $U_2, \cdots, U_n$ using the formulas in \eqref{e514}. However the formulas for $U_k$, $2\leq k\leq n$ are quite complicated in general. 

\medskip
In this appendix, we focus on the case of $SU(3)$ and give   the explicit formulas for $n = 2$. Consider 
\begin{equation}
\label{toda2s}
-\Delta u_1 = 2e^{u_1} - e^{u_2} - 4\pi\gamma_1\delta_0, \;\; -\Delta u_2 = 2e^{u_2} - e^{u_1} - 4\pi\gamma_2\delta_0 \quad \text{in }\mathbb R^2, \quad  \int_{\mathbb R^2}e^{u_i} <\infty, \;\; i = 1, 2,
\end{equation}
with $\gamma_1, \gamma_2 > -1$. Our result is
\BT
\label{main2}
Assume that $(u_1, u_2)$ is solution of \eqref{toda2s}. 
\begin{itemize}
\item If $\gamma_1, \gamma_2 \in \N$. The solution space is an eight dimensional smooth manifold. More precisely, we have
\begin{align}
\label{classSU}
e^{u_1} = 4\Gamma|z|^{2\gamma_1}\frac{Q}{P^2}, \quad e^{u_2} = 4\Gamma|z|^{2\gamma_2}\frac{P}{Q^2} \quad \mbox{in } \C
\end{align}
with $\Gamma = (\gamma_1 + 1)(\gamma_2 + 1)(\gamma_1 + \gamma_2 + 2)$ and
\begin{align*}
P(z) = & \; (\gamma_2 + 1)\xi_1 + (\gamma_1 + \gamma_2 + 2)\xi_2\left|z^{\gamma_1 + 1} - c_1\right|^2 +
\frac{\gamma_1 +1}{\xi_1\xi_2}\left|z^{\gamma_1 + \gamma_2 + 2} - c_2z^{\gamma_1 + 1} - c_3 \right|^2,\\
Q(z) = & \; (\gamma_1 + 1)\xi_1\xi_2 + \frac{\gamma_1 + \gamma_2 + 2}{\xi_2}\left|z^{\gamma_2 + 1} - \frac{(\gamma_1 + 1)c_2}{\gamma_1 + \gamma_2 + 2}\right|^2\\
& \;+ \frac{\gamma_2 + 1}{\xi_1}\left|z^{\gamma_1 + \gamma_2 + 2} - \frac{(\gamma_1 + \gamma_2 + 2)c_1}{\gamma_2 + 1}z^{\gamma_1 + 1} + \frac{(\gamma_1 + 1)c_3}{\gamma_2 + 1} \right|^2,
\end{align*}
where $c_1, c_2, c_3 \in \C$, $\xi_1, \xi_2 > 0$.
\item If now $\gamma_1 \not \in \N$, $\gamma_2 \not \in \N$ and $\gamma_1+\gamma_2 \not \in \mathbb Z$, then $c_1 = c_2 = c_3 = 0$, the solution manifold to \eqref{toda2s} is of two dimensions.
\item If $\gamma_1 \in \N$, $\gamma_2 \not \in \N$, then $c_2 = c_3 = 0$; if $\gamma_1 \not\in \N$, $\gamma_2 \in \N$, there holds $c_1 = c_3 = 0$; we get a four dimensional solution manifold in both cases.
\item If $\gamma_1 \not\in \N$, $\gamma_2 \not \in \N$ but $\gamma_1 + \gamma_2 \in \mathbb Z$, then $c_1 = c_2 = 0$, the solution manifold to \eqref{toda2s} is of four dimensions.
\end{itemize}
In all cases, we have 
\begin{equation}
\label{quan}
\int_{\R^2} e^{u_1} dx = \int_{\R^2} e^{u_2} dx =4 \pi (\gamma_1+\gamma_2 +2).
\end{equation}
\ET

\medskip
The proof can follow directly from the formulas \eqref{e107n} and \eqref{e514}. Here in the below we give direct calculations instead of the general consideration in section 5. 

\medskip
Define $(U_1, U_2)$ and $\alpha_1, \alpha_2$ by \eqref{singSU}. Denoting
$$W_1 = - e^{U_1}\left(e^{-U_1}\right)''' = U_{1, zzz} - 3U_{1, zz}U_{1, z} +
U_{1, z}^3,$$
then $W_{1, \bar z} = -U_{1, z\bar z}\left[U_{1, zz} + U_{2, zz} - U_{1,
z}^2 - U_{2, z}^2 + U_{1, z}U_{2, z}\right] :=  -U_{1, z\bar z}W_2$. As before, we can claim that $W_{2, \bar z} = 0$ in $\C^*$. By stuyding the behavior of $W_2$ at $\infty$, we get
\begin{align*}
W_2 = \frac{w_2}{z^2} \;\; \mbox{in $\C^*$}\quad \mbox{where } \; w_2 = -\alpha_1^2 - \alpha_2^2 + \alpha_1\alpha_2 -\alpha_1 - \alpha_2.
\end{align*}
As $(W_1 + U_{1, z}W_2)_{\bar z} = U_{1, z}W_{2, \bar z} = 0$ in $\C^*$, by considering $z^3(W_1 + U_{1, z}W_2)$, there holds
$$W_1 + U_{1, z}W_2 = \frac{w_1}{z^3}  \;\; \mbox{in $\C^*$}\quad \mbox{where } \; w_1 = 2\alpha_1 + 3\alpha_1^2 + \alpha_1^3 + \alpha_1w_2.$$
Combine these informations, the function $f := e^{-U_1}$
satisfies
\begin{align}
\label{fSU3}
f_{zzz} = -fW_1 = -\frac{w_1}{z^3}f + fU_{1, z}\frac{w_2}{z^2} = -\frac{w_2}{z^2}f_z - \frac{w_1}{z^3}f \quad \mbox{in } \C^*.
\end{align}
Consider special solution of \eqref{fSU3} like $z^\beta$, then $\beta$ should satisfy $\beta(\beta - 1)(\beta - 2) + w_2\beta + w_1 = 0$. We check readily that the equation of $\beta$ has three roots: $\beta_1 = -\alpha_1$, $\beta_2 = \alpha_1 + 1 - \alpha_2$ and $ \beta_3 = \alpha_2 + 2$. Hence $\beta_3 - \beta_2 = \gamma_2 + 1 > 0$ and $\beta_2 - \beta_1 = \gamma_1 + 1 > 0$.  We obtain finally $f(z)  = \sum_{1\leq i, j \leq 3} b_{ij}z^{\beta_i}\bar z^{\beta_j}$ with an hermitian matrix $(b_{ij})$. 

\medskip
In the following, we show how to get explicit formulas of $U_i$ for just two cases, and all the others can be treated similarly. The formulas of $u_i$ or the quantization \eqref{quan} of the integrals are clearly direct consequences of the expressions of $U_i$.

\begin{itemize}
 \item {\it Case 1}: $\gamma_i \notin \N$ and $\gamma_1 + \gamma_2 \notin \mathbb Z$.
\end{itemize}
To get a well defined real function $f$ in $\C^*$, we have $b_{ij} = 0$ for $i\ne j$, so that
$$f = e^{-U_1} = \sum_{i = 1}^3 a_i  |z|^{2\beta_i} \;\; \mbox{in } \C^*, \quad \mbox{with } a_i \in \R.$$
Therefore direct calculation yields
\begin{align*}
\frac{e^{-U_2}}{4} = -e^{-2U_1}U_{1, z\bar z} = ff_{z\bar z} - f_zf_{\bar z} = \sum_{1 \leq i < j \leq
3}a_ia_j(\beta_i - \beta_j)^2|z|^{2(\beta_i + \beta_j - 1)}.
\end{align*}
Moreover, there holds also $e^{-U_1} =
-4e^{-2U_2}U_{2, z\bar z}$. With the explicit values of
$\beta_i$, we can check that $(U_1, U_2)$ is a solution if and only if
\begin{align}
\label{ai} a_1a_2a_3\Gamma^2 = \frac{1}{64}\;\; \mbox{where $\Gamma= (\gamma_1 + 1)(\gamma_2 + 1)(\gamma_1 + \gamma_2 + 2)$} ,
\end{align}
or equivalently
$$a_1 = \frac{(\gamma_2 +1)\xi_1}{4\Gamma}, \quad a_2 = \frac{(\gamma_1 + \gamma_2 + 2)\xi_1}{4\Gamma}, \quad a_3 = \frac{(\gamma_1 +1)}{4\Gamma\xi_1\xi_2} \quad \mbox{with }\; \xi_1, \xi_2 > 0.$$
Indeed, the positivity of $e^{-U_1}$ in $\C^*$ implies that $a_1, a_3 > 0$, so is $a_2$ by \eqref{ai}.

\begin{itemize}
 \item {\it Case 2}: $\gamma_1 \in \N$ but $\gamma_2 \not\in \N$. 
\end{itemize}
We get then
$$e^{-U_1} = \sum_{i = 1}^3 a_i  |z|^{2\beta_i} + \frac{{\rm Re}\left(\lambda
z^{\gamma_1 + 1}\right)}{|z|^{2\alpha_1}}\;\; \mbox{in } \C^*, \quad \mbox{with } a_i \in \R,\; \lambda \in \C.$$
If $a_2 \ne 0$, changing
eventually the value of $a_1$, there exists $c_1 \in \C$ such that
$$e^{-U_1} = \frac{a_1 + a_2\left|z^{\gamma_1 + 1} - c_1\right|^2 + a_3 |z|^{2(\gamma_1 + \gamma_2 +
2)}}{|z|^{2\alpha_1}} \quad \mbox{in } \C^*.$$ We obtain then the expression of
$e^{-U_2}$ directly and we can check that the necessary and sufficient condition required to get solutions of \eqref{toda2s} is always \eqref{ai}. We leave the details for interested readers. This yields
$$e^{-U_1} = \frac{1}{4\Gamma|z|^{2\alpha_1}}\left[(\gamma_2 + 1)\xi_1 + (\gamma_1 + \gamma_2 + 2)\xi_2\left|z^{\gamma_1 + 1} - c_1\right|^2 +
\frac{\gamma_1 +1}{\xi_1\xi_2}|z|^{2(\gamma_1 + \gamma_2 +
2)}\right]$$ and \begin{align*} e^{-U_2} = &
\;\;\frac{1}{4\Gamma|z|^{2\alpha_2}}\left[(\gamma_1 +
1)\xi_1\xi_2 + \frac{\gamma_1 + \gamma_2 + 2}{\xi_2}|z|^{2(\gamma_2 +
1)} + \frac{\gamma_2 +1}{\xi_1}|z|^{2(\gamma_2 + 1)}\left|z^{\gamma_1
+ 1} - \frac{(\gamma_1 + \gamma_2 + 2)c_1}{\gamma_2 +
1}\right|^2\right]. \end{align*} 

So it remains to eliminate the case $a_2 = 0$. If $a_2 = 0$, we can rewrite $$f = \frac{a_1 + {\rm Re}\left(\lambda
z^{\gamma_1 + 1}\right) + a_3|z|^{2(\gamma_1 + \gamma_2 +
2)}}{|z|^{2\alpha_1}}\;\; \mbox{in } \C^*$$
where $\lambda \in \C$. Direct calculation yields
\begin{align*}
\frac{e^{-U_2}}{4} = ff_{z\bar z} - f_zf_{\bar z} = |z|^{2(-\alpha_2 + \gamma_2 + 1)}\left[c_1'|z|^{-2(\gamma_2+1)} + c_2' + c_3'{\rm Re}\left(\lambda
z^{\gamma_1 + 1}\right)\right]
\end{align*}
where
\begin{align*}
c_1' = -\frac{|\lambda|^2(\gamma_1 + 1)^2}{4}, \quad c_2' = a_1a_3(\gamma_1 + \gamma_2 +
2)^2, \quad c_3' = a_3(\gamma_1 + \gamma_2 +
2)(\gamma_2 + 1).
\end{align*}
As $e^{-U_2} > 0$, we must have $c_1' \geq 0$. So we get $\lambda = 0$, and we find the expression of $f$ as in {\sl Case 1} with $a_2 = 0$. Then we need to verify the equation \eqref{ai}. However this is impossible since $a_2 = 0$. Thus $a_2$ must be nonzero.

\end{document}